\documentclass[12pt]{article}
\usepackage[margin=0.75in,top=0.6in]{geometry}
\usepackage[utf8]{inputenc}
\usepackage{amsmath,amsfonts,amssymb}
\usepackage{algorithm}
\usepackage{algpseudocode}
\usepackage{graphicx}
\usepackage{subfig}
\usepackage{float}
\usepackage{multirow} 
\usepackage{array}  
\usepackage{booktabs}  
\usepackage{siunitx}  
\usepackage{caption}
\usepackage{hyperref}
\usepackage{nomencl}
\usepackage{multicol}
\makenomenclature
\renewcommand{\nompreamble}{\begin{multicols}{2}}
\renewcommand{\nompostamble}{\end{multicols}}
\date{}
\title{Multi-Scenario and Stochastic Thermo-Electro-Mechanical Modelling of Failure in Power Transmission Lines}
\begin{document}
\author{
    Prakash KC\thanks{Department of Mechanical Engineering, Michigan State University} \and
    Maryam Naghibolhosseini\thanks{Department of Communicative Sciences and Disorders, Michigan State University} \and
    Mohsen Zayernouri\thanks{Department of Mechanical Engineering \& Department of Statistics and Probability, Michigan State University, Corresponding Author; zayern@msu.edu}
}
\captionsetup{font=small}
\maketitle
\begin{abstract}
    Transmission lines, crucial to the power grid, are subjected to diverse environmental conditions such as wind, temperature, humidity, and pollution. While these conditions represent a consistent impact on the transmission lines, certain unpredictable conditions such as unexpected high wind, wildfire, and icing pose catastrophic risks to the reliability and integrity of the transmission lines. These factors in the presence of initial damage and electrical loads greatly affect the material properties. In this paper, we develop a comprehensive thermo-electro-mechanical model to investigate the long-term effect of unexpected high wind, wildfire, and ice on transmission lines. This study offers an in-depth perspective on temperature and damage evolution within the power lines by incorporating a phase field model for damage and fatigue, alongside thermal and electrical models. We define a state function to assess the failure, considering damage and temperature. We study three scenarios deterministically to establish a basic understanding and analyze the stochastic behavior using the Probabilistic Collocation Method (PCM). We utilize PCM for forward uncertainty quantification, conducting sensitivity analysis, and evaluating the probability of failure. This approach offers an in-depth examination of the potential risks associated with transmission lines under unfavorable circumstances.

\end{abstract}   
{
\hspace{10mm}\noindent \textbf{Keywords:}
\footnotesize{Transmission line, Probability of Failure, Finite Element Method, Probabilistic Collocation Method, Sensitivity Analysis, Wildfire, Icing}\\
}
\nomenclature{\(\theta_{lim}\)}{Threshold Temperature}
\nomenclature{\(\theta_{max}\)}{Maximum Temperature}
\nomenclature{\(\varphi_{lim}\)}{Threshold Damage}
\nomenclature{\(\varphi_{max}\)}{Maximum Damage}
\nomenclature{\(A_{0}\)}{Mean Coefficient}
\nomenclature{\(A_{n}, B_{n}\)}{Frequency Coefficients}
\nomenclature{\(t\)}{Time}
\nomenclature{\(\mathcal{F}\)}{Fatigue Field}
\nomenclature{\(u\)}{Displacement}
\nomenclature{\(Y\)}{Young Modulus}
\nomenclature{\(g_c\)}{Fracture Energy Release Rate}
\nomenclature{\(\gamma\)}{Phase Field Layer Width}
\nomenclature{\(\mathcal{H}'(\varphi), \mathcal{H}_f'(\varphi) \)}{Damage Potentials}
\nomenclature{\(\theta_{c}\)}{Conductor Temperature}
\nomenclature{\(\varphi\)}{Damage}
\nomenclature{\(\theta_{c}\)}{Conductor Temperature}
\nomenclature{\(\theta_{0}\)}{Reference Temperature}
\nomenclature{\(\rho\)}{Material Density}
\nomenclature{\(a\)}{Aging Rate}
\nomenclature{\(q_j\)}{Joule Heating}
\nomenclature{\(\varphi\)}{Damage}
\nomenclature{\(q_c\)}{Convection Cooling}
\nomenclature{\(h\)}{Convective Heat Transfer Coefficient}
\nomenclature{\(v\)}{Wind Velocity}
\nomenclature{\(\rho_{air}\)}{Air Density}
\nomenclature{\(\rho_{ice}\)}{Ice Density}
\nomenclature{\(\nu\)}{Kinematic Viscosity of Air}
\nomenclature{\(D\)}{Conductor Diameter}
\nomenclature{\(Re_{D}\)}{Reynolds Number}
\nomenclature{\(Pr\)}{Prandtl Number}
\nomenclature{\(I_b\)}{Base Current}
\nomenclature{\(\theta_b\)}{Base Temperature}
\nomenclature{\(w_b\)}{Base Wind Speed}
\nomenclature{\(q_j\)}{Joule Heating}
\nomenclature{\(q_f\)}{Radiative Heat Transfer}
\nomenclature{\(q_i\)}{Heat Transfer due to ice}
\nomenclature{\(\sigma_e\)}{Electrical Conductivity}
\nomenclature{\(\sigma_{E,T}\)}{Non Degraded Electrical Conductivity}
\nomenclature{\(\sigma_0\)}{Electrical Conductivity at Reference Temperature}
\nomenclature{\(P_w\)}{Wind Pressure}
\nomenclature{\(w_w\)}{Wind Load}
\nomenclature{\(C_D\)}{Drag Coefficient}
\nomenclature{\(\alpha\)}{Span Factor}
\nomenclature{\(\theta_w\)}{Wind Direction}
\nomenclature{\(w\)}{Test Function}
\nomenclature{\(A(x)\)}{Cross section Area}
\nomenclature{\(A_0\)}{Undamaged Cross section Area}
\nomenclature{\(J\)}{Current Density}
\nomenclature{\(A_\sigma\)}{Level of Damage}
\nomenclature{\(\delta_t\)}{Time step}
\nomenclature{\(h_B\)}{Bernoulli Random Variable}
\nomenclature{\(p_h\)}{Bernoulli parameter}
\nomenclature{\(H_0\)}{Initial Horizontal Tension}
\nomenclature{\(p_f\)}{Probability of Failure}
\nomenclature{\(\sigma\)}{Stefan Boltzmann Constant}
\nomenclature{\(\epsilon\)}{Flame Emissivity}
\nomenclature{\(\theta_f\)}{Radiating Surface average temperature}
\nomenclature{\(\phi\)}{view factor} 
\nomenclature{\(\kappa\)}{Thermal Conductivity of Condutor} 
\nomenclature{\(\kappa_{air}\)}{Thermal Conductivity of Air}
\nomenclature{\(\kappa_{ice}\)}{Thermal Conductivity of Ice}
\nomenclature{\(\theta_{ice}\)}{Temperature of Ice}
\nomenclature{\(r_1\)}{Inner Radius of Ice Layer}
\nomenclature{\(r_2\)}{Outer Radius of  Ice Layer}
\nomenclature{\(m_e\)}{Mass of Melted Ice}
\nomenclature{\(L_e\)}{Latent Heat of Fusion}
\nomenclature{\(w_c\)}{Conductor weight per unit length}
\nomenclature{\(w_i\)}{Ice weight per unit length}
\nomenclature{\(V\)}{Volume of Ice per unit length}
\nomenclature{\(t_{ice}\)}{Ice Thickness}
\nomenclature{\(g\)}{Acceleration due to gravity}

\printnomenclature

\section{Introduction}
The power grid is a complex, highly interconnected network due to which the failure of even a single component can lead to cascading failures \cite{albert2004structural,crucitti2004model,kinney2005modeling}. Natural phenomena such as storms, hurricanes, tornadoes, wildfires, icing, and extreme temperatures are the primary factors of these failures \cite{ward2013effect,gao2018potential}. Given the complexity and evident sensitivity of the power grid, accurately predicting system failures is crucial. 
It is difficult yet important to identify the sensitive areas and analyze failure, especially considering the significance of these areas for energy distribution and national security. Power grids in such sensitive areas are influenced by a variety of unpredictable factors like weather, damage, and aging.
Transmission lines of the power grid are the most vulnerable component due to their exposure to dynamic weather conditions.  Detailed studies have been made on the reliability of transmission lines. The effect of thermal stress on the life span of transmission lines, specific to Aluminum Conductor Steel-Reinforced (ACSR) was investigated in \cite{hathout2018impact}. In studies such as \cite{alminhana2023transmission,dua2015dynamic,stengel2014finite,hamada2011behaviour}, the dynamic effect of wind was investigated. Additionally, the combined effect of ice and wind loads on transmission lines was examined in recent studies \cite{yang2013probability, rossi2020combined}. Radiative heat due to wildfires, contributing to thermal failure through overheating, was also addressed in specific studies \cite{koufakis2010wildfire,guo2018determination}.
Furthermore, the influence of ambient temperature, wind speed, and electrical current on conductor temperatures has been thoroughly investigated \cite{bockarjova2007transmission}. Other research has focused on how various weather conditions affect the thermal ratings of conductors \cite{bendik2018influence,pytlak2009precipitation,castro2017study}.

All these studies are crucial for understanding the failure mechanisms in power transmission lines. However, they often have limitations due to the narrow focus on specific weather conditions, relying on thermal and electrical models, or are typically based on short-term simulations. There is a clear need for developing models that integrate the electrical and thermal properties of overhead transmission lines with advanced damage and fatigue mechanisms. Moreover, the model must take into account the effects of catastrophic events such as unexpected seasonal high winds, wildfires, and icing conditions. Such models could significantly enhance the accuracy of predicting transmission line failures over the long term.

Recently, phase-field models have emerged as a solid area of research in damage and fatigue modeling. These models address a range of issues from brittle \cite{miehe_phase_2010,miehe_thermodynamically_2010,borden_higher-order_2014,de2021data} and ductile fractures \cite{ambati_phase-field_2015,ambati_phase-field_2016} to dynamic \cite{borden_phase-field_2012,hofacker_phase_2013} and non-isothermal fatigue fractures \cite{boldrini_non-isothermal_2016}. The phase-field model effectively handles crack initiation and growth through a smooth crack representation, eliminating the need for explicit crack geometry tracking. Recent studies have integrated phase-field models with electrical and thermal simulations to investigate phenomena such as polymer breakdown under alternating voltage \cite{zhu2020phase}, transitions in the state from paraelectric to ferroelectric \cite{woldman2019thermo}, and the failure of polymer-based dielectrics \cite{shen2019phase}. However, the application of these models remains limited. Other continuum damage models have been applied to electrical conductors \cite{kaiser2021fundamentals}, self-sensing materials \cite{nayak2019microstructure}, and the thermo-electro-mechanical degradation of electrical contacts \cite{shen2021numerical}. However, there is a notable gap in modeling the impact of damage and fatigue on accurately predicting the life span of transmission lines. Additionally, integrating stochastic analysis to account for variables like weather conditions and existing damage could significantly improve the accuracy of these predictions.

The failure analysis of the transmission line requires an understanding of the interconnection between mechanical, thermal, and electrical components, along with the inherent uncertainties. The phase field model is widely accepted for damage evolution modeling, while dislocation dynamics can capture the details of the microstructural mechanisms \cite{de2021atomistic}. Studies by Chhetri et al. \cite{chhetri2023comparative} and De et al. \cite{de2023machine} have shown that material behavior can be significantly affected by dislocation interactions, which result in failure under dynamic conditions.

Under high temperatures, material behavior deviates from normal to more complex behavior, including both viscous and elastic components. Fractional visco-elasto-plastic models \cite{suzuki2022general} can effectively capture these behaviors, leading to better structural analysis \cite{suzuki2016fractional}, accurate damage evolution \cite{suzuki2021thermodynamically}, and improved large-scale behavior \cite{suzuki2023fractional}. Integrating these fractional-order models can enhance the failure prediction of transmission lines. Additionally, incorporating the fractional order as a random variable can better capture the system response \cite{kharazmi2019operator}.

The traditional Monte Carlo method, a benchmark in stochastic analysis for computing the Quantity of Interest (QoI), is known for its simplicity \cite{fishman2013monte,smith2013uncertainty}. However, it has a slow convergence rate and requires numerous realizations making the process computationally intensive. This highlights the need for more efficient computational approaches in stochastic analysis. Although there are alternative established methods, such as Polynomial Chaos \cite{xiu2002modeling,xiu2002wiener} and its extension through Galerkin projection \cite{stefanou2009stochastic,babuvska2005solving,babuska2004galerkin}, they often require changes to the original equations, making them intrusive and less practical for complex scenarios. In that case, the non-intrusive Probabilistic Collocation Method (PCM) \cite{babuvska2007stochastic,xiu2005high} can be an excellent choice. PCM maintains the original solution structure and allows for independent sampling of realizations, thus offering better convergence rates than the Monte Carlo method. Despite the challenge of dimensionality, PCM's effectiveness can be significantly improved with methods such as sparse grids \cite{smolyak1963quadrature} and active subspace techniques \cite{constantine2017global,constantine2015exploiting,constantine2014active}.
The Probabilistic Collocation Method (PCM) has been effectively implemented to explore the Uncertainty quantification in power grid systems \cite{hockenberry2004evaluation,lin2014uncertainty}, however, these studies primarily address the short-term electrical behaviors. In the area of damage phase-field models, PCM has proven effective for conducting forward uncertainty quantification (UQ) and sensitivity analysis (SA) \cite{barros2021integrated}.

In this study, we introduce a coupled thermo-electro-mechanical model incorporating equations for displacements, material damage, fatigue, temperature, and voltage similar to our previous work in \cite{demoraes2024thermoelectromechanicalmodellongtermreliability,kc2024thermoelectromechanicalmodelingpowertransmission}. We integrate wind, temperature, and current demand in our model to analyze the failure of the transmission line in the presence of initial damage. We then incorporate unexpected weather conditions such as high wind, wildfire, and icing on our model to understand their impact on the reliability of transmission lines. Instead of using probability distributions for failure parameters, we apply the PCM to propagate parametric uncertainty to the temperature or damage outputs of our physics-based material model. We assess the impact of each stochastic parameter using variance-based global sensitivity analysis and use PCM to estimate the probability of failure over time.
The key contributions of this research are:

\begin{itemize}
    \item We introduce an integrated model that combines phase-field modeling for damage and fatigue with thermal and electrical effects.
    \item We extend the use of the Probabilistic Collocation Method (PCM) beyond standard uncertainty and sensitivity analyses, applying it to directly calculate the probability of failure. This is achieved by transforming the limit state function into a failure-indicator Bernoulli random variable.
    \item We explore the long-term behavior of the coupled physical system, identifying influential factors that lead to the early failure of transmission lines under various scenarios, from unexpected high-season wind to wildfire and icing.
\end{itemize}
The structure of this paper is as follows: Section 2. introduces the Problem Statement and elaborates the several representative scenarios. In Section 3, we discuss the thermo-electro-mechanical model for transmission line failures elaborating on each model in detail. In section 4, we detail the deterministic analysis using the Finite-Element Method. Section 5 covers the stochastic methods, presenting the Probabilistic Collocation Method (PCM) for uncertainty, sensitivity, and probability of failure analyses. Concluding remarks are provided in Section 5.

\section{Problem Statement}
\label{sec:Problem Statement}

Operating temperature and damage are important indicators for the safe operation of transmission lines. High temperatures may lead to annealing and sagging of transmission wires, while cyclic loading, especially on already damaged wires, can result in cracks initiation and propagation, ultimately leading to failure in transmission lines. Often, these issues go unnoticed until a cable failure occurs. Environmental factors such as ambient air temperature, wind speed, and increased current demand in the presence of initial damage can negatively influence these parameters. Furthermore, the unexpected seasonal wind, wildfire, and icing pose adverse effects on these parameters, pushing them beyond threshold limits, and ultimately resulting in early failure.

The operating temperature and damage of the conductor are, therefore the primary focus of our analysis. The interaction among various physical effects, loading conditions, material parameters, and unexpected environmental factors will determine the long-term failure of the transmission line. We aim to explore the continuous effect of multiphysics on the operating temperature and damage, particularly in the presence of initial damage. 

\subsection{Representative Scenarios}
We consider three different scenarios for the detailed study of the effect of unexpected season wind, wildfire, and icing. For our analysis, we use wind and temperature data for Texas, California, and Alaska from the \cite{noaa_ccd_2024} and \cite{nws_2024} to study the respective scenarios. For simplification, we consider monthly average data of wind and temperature and use the Discrete Fourier Transform (DFT) and Fourier analysis to obtain continuous data for a year and use it as loading conditions for each year in our study. The DFT transform a sequence \(x_n\) into:
\[X_k = \sum_{n=0}^{N-1} x_n e^{-i2\pi \frac{kn}{N}}\]
where \(x_n\) is the nth data sample, and \(N\) is the total number of samples. 

The component \(A_0\) is the mean of the data, and coefficients \(A_n\) and \(B_n\) for \(n=1, 2, ..., N/2\) are the frequencies obtained from the real and imaginary parts of \(X_k\). The cyclic loading equation is reconstructed using the Fourier series:
\[f(t) = A_0 + \sum_{n=1}^{N/2} [A_n \cos(2\pi n \frac{t}{T}) + B_n \sin(2\pi n \frac{t}{T})]\]
where $t$ is a scaled time matching the original time domain, and $T$ is the period of the data set. The time is scaled in such a way that 100 iterations correspond to the cyclic loading condition of 12 months.
\begin{itemize}
	\item \textbf{Scenario 1 - High Seasonal Winds and Temperature:} To study the effect of unexpected high wind in the regions with high variation of temperature within a year, we consider the wind and temperature of Texas (Amarillo). Initially, we analyze the failure without the effect of unexpected high wind in the presence of different levels of initial damage. Later, we introduce unexpected seasonal high winds to study the effect on the failure of an initially damaged transmission line. For this scenario, high wind speeds of 50, 75, and 100 (ft/s) are considered. The original data of wind and temperature are given in Table.\ref{tab:Texas data}. The original and transformed data for wind and temperature are shown in Fig.\ref{fig:First}.
   \begin{table}[ht]
   \caption{Original wind \cite{noaa_ccd_2024} and temperature \cite{nws_2024} data for Amarillo, TX }
   \small
    \begin{tabular}{cc*{12}{p{0.8cm}}}
        \hline
        Month & Jan & Feb & Mar & Apr & May & Jun & Jul & Aug & Sep & Oct & Nov & Dec \\
        \hline
        Wind Speed (ft/s) & 17.75 & 18.92 & 20.39 & 21.56 & 20.09 & 20.39 & 18.19 & 16.72 & 17.75 & 18.33 & 18.63 & 17.89 \\
        Temperature (K) & 276.87 & 280.98 & 284.76 & 287.32 & 290.43 & 298.32 & 301.37 & 297.59 & 295.15 & 291.87 & 284.37 & 276.54 \\
        \hline
    \end{tabular}
    \label{tab:Texas data}
\end{table}

 \begin{figure}[H]
    \centering
    % Row for wind dataset
    {\includegraphics[height=6cm,width=0.45\textwidth]{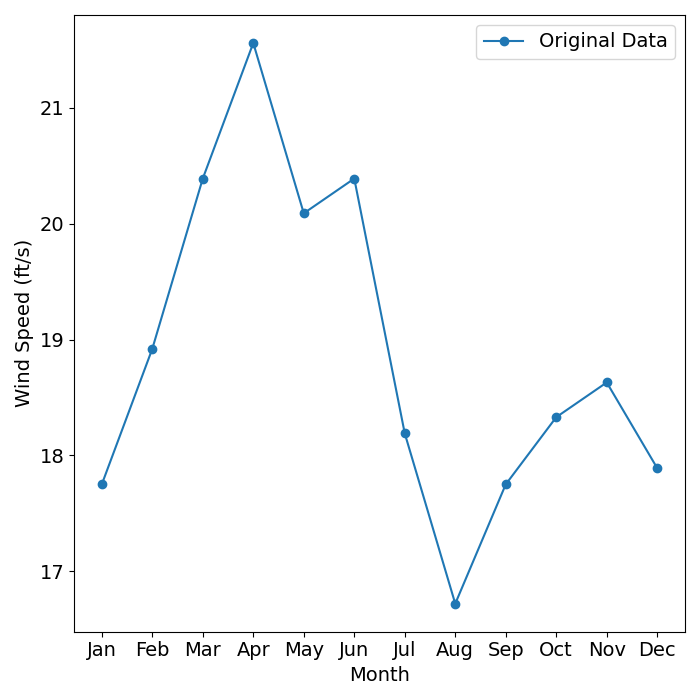}}
    {\includegraphics[height=6cm,width=0.45\textwidth]{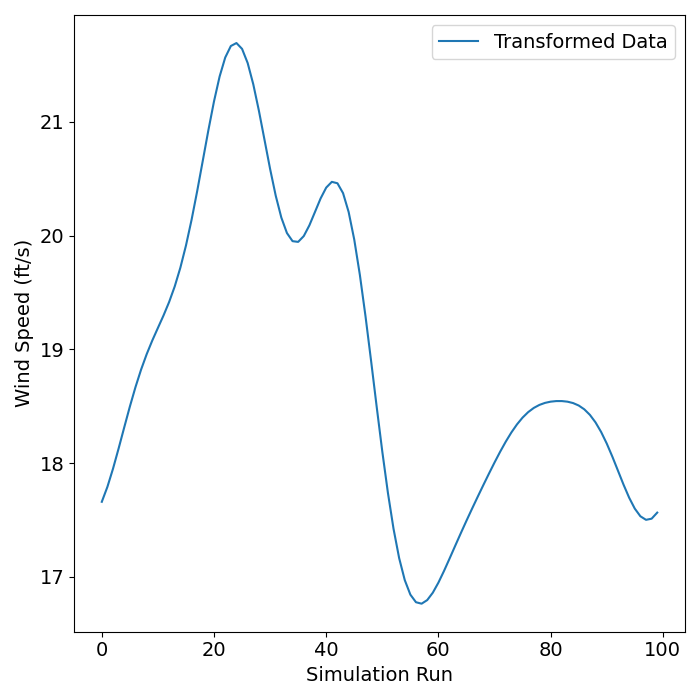}}
   \hfill
    {\includegraphics[height=6cm,width=0.45\textwidth]{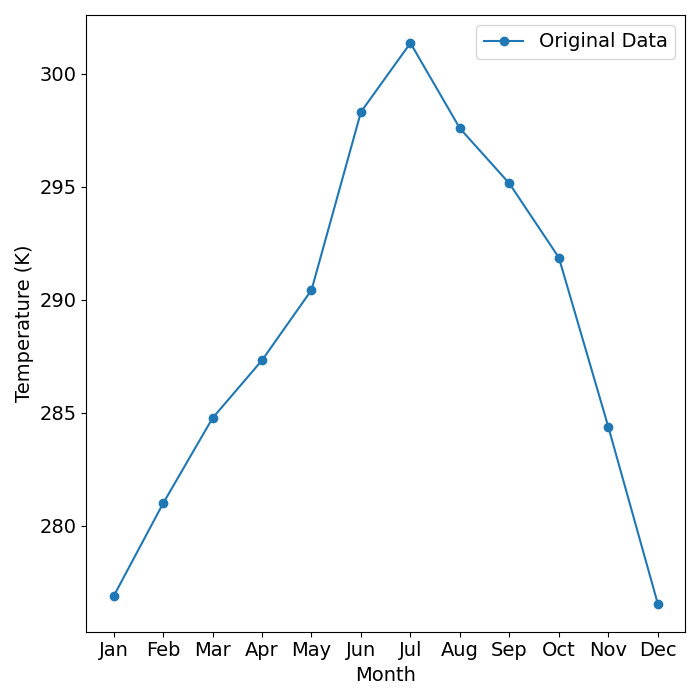}}
    {\includegraphics[height=6cm,width=0.45\textwidth]{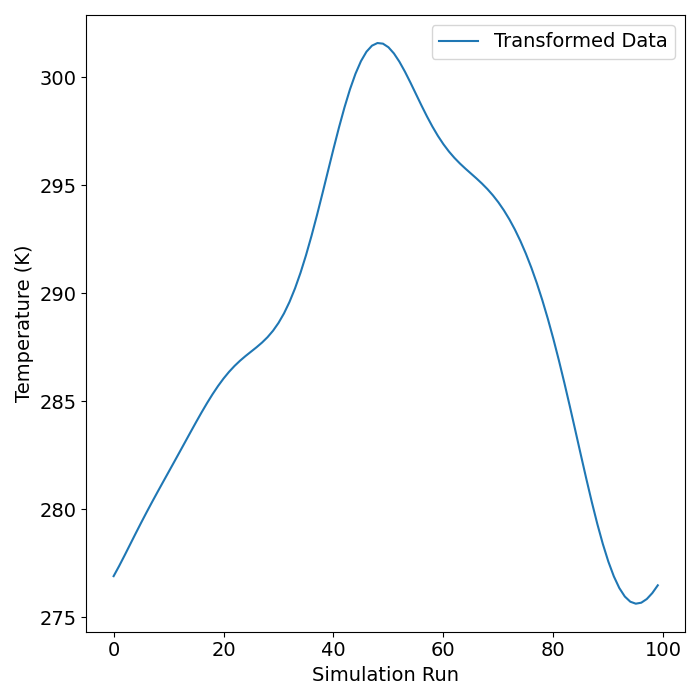}}
    \caption{Discrete and continuous variation of wind and temperature for Amarillo, TX}
    \label{fig:First}
   \end{figure} 
 
	\item \textbf{Scenario 2 - Wildfire:} In this scenario, we consider the effect of wildfire on the life span of the transmission line. According to the \href{https://www.energy.ca.gov/sites/default/files/2019-11/Energy_CCCA4-CEC-2018-002_ADA.pdf}{(California Energy Commission, 2018)}, due to the effect of wildfire, the state incurred costs of around \$700 million in transmission and distribution-related damages during the period from 2000 to 2016. Therefore, we study the effect of wildfire on the transmission lines using the wind and temperature data of California (San Diego). We first analyze the failure of the transmission line without the effect of wildfire in the presence of initial damage. Then, we consider the effect of wildfire on the initially damaged transmission line varying the distance between the transmission line and the wildfire in the form of a view-factor. Table.\ref{tab:California data} shows the wind and temperature data of California. The original and transformed data for wind and temperature are shown in Fig.\ref{fig:Second}.
 \begin{table}[ht]
 \caption{Original wind \cite{noaa_ccd_2024} and temperature \cite{nws_2024} data for San Diego, CA}
   \small
    \begin{tabular}{cc*{12}{p{0.8cm}}}
        \hline
        Month & Jan & Feb & Mar & Apr & May & Jun & Jul & Aug & Sep & Oct & Nov & Dec \\
        \hline
        Wind Speed (ft/s) & 7.48 & 8.95 & 9.68 & 10.56 & 10.85 & 10.56 & 10.27 & 9.97 & 9.68 & 8.36 & 7.48 & 7.19 \\
        Temperature (K) & 289.15 & 290.43 & 292.37 & 291.93 & 290.93 & 293.65 & 295.59 & 297.09 & 298.26 & 296.71 & 290.54 & 287.65 \\
        \hline
    \end{tabular}
    \label{tab:California data}
\end{table}

 \begin{figure}[H]
    \centering
    % Row for wind dataset
    {\includegraphics[height=6cm,width=0.45\textwidth]{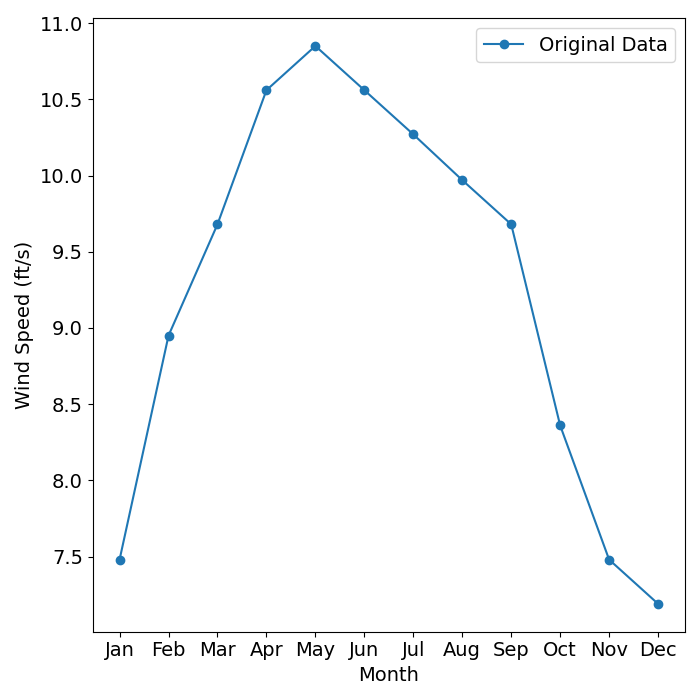}}
    {\includegraphics[height=6cm,width=0.45\textwidth]{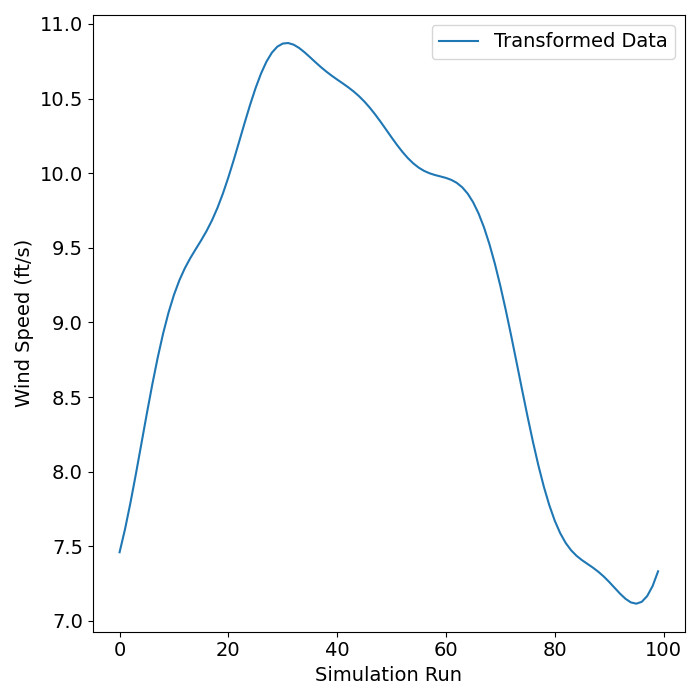}}
   \hfill
    {\includegraphics[height=6cm,width=0.45\textwidth]{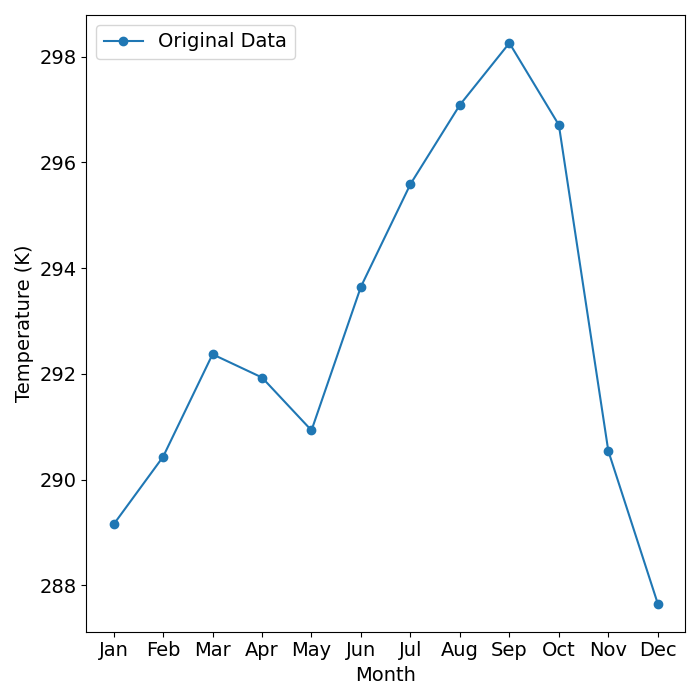}}
    {\includegraphics[height=6cm,width=0.45\textwidth]{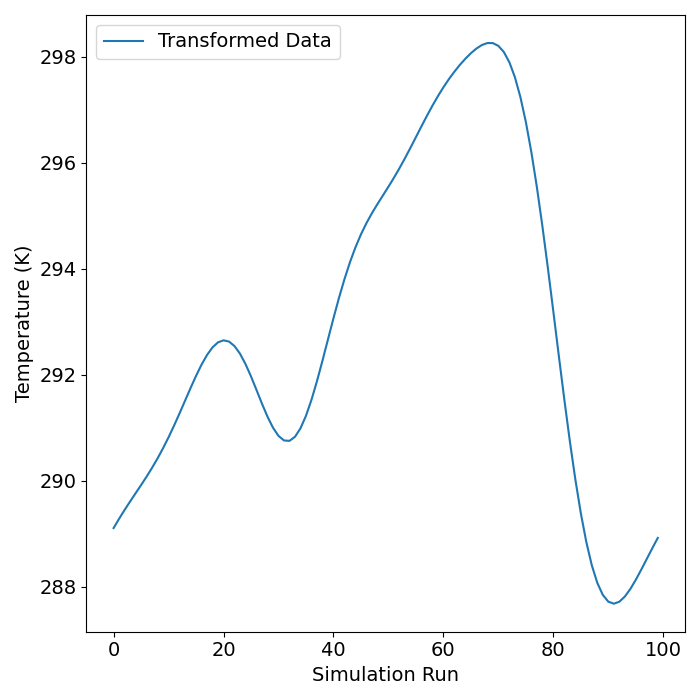}}
    \caption{Discrete and continuous variation of wind and temperature for San Diego, CA}
    \label{fig:Second}
   \end{figure} 
	\item \textbf{Scenario 3 - Icing:} In the third scenario, we study the effect of a uniform ice layer around the transmission line. The ice layer on the transmission lines can have several effects, however, here, we focus on the effect of ice in terms of the mechanical load on the transmission line, heat transfer, and the electrical conductivity of the wire. To study the effect of icing on the life span of a transmission line, we consider the wind and temperature data of Alaska (Bethel), which is one of the coldest and windiest states. We first study the failure on the transmission lines in the presence of different levels of initial damage without considering the effect of ice. Then, we assume the uniform layer of ice with the thickness of 0.1, 0.25, and 0.5 inch around the cables according to the design criteria of \href{https://www.rd.usda.gov/files/UEP_Bulletin_1724E-200.pdf}{NESC}. Table.\ref{tab:Alaska data} shows the wind and temperature data of Alaska. The original and transformed data for wind and temperature are shown in Fig.\ref{fig:Third}.
 \begin{table}[ht]
 \caption{Original wind \cite{noaa_ccd_2024} and temperature \cite{nws_2024} data for Bethel, AK}
   \small
    \begin{tabular}{cc*{12}{p{0.8cm}}}
        \hline
        Month & Jan & Feb & Mar & Apr & May & Jun & Jul & Aug & Sep & Oct & Nov & Dec \\
        \hline
        Wind Speed (ft/s) & 20.68 & 20.24 & 18.77 & 17.01 & 15.55 & 14.23 & 14.08 & 15.11 & 15.40 & 16.72 & 18.04 & 18.92 \\
        Temperature (K) & 256.48 & 261.82 & 257.32 & 268.98 & 279.54 & 283.65 & 284.98 & 284.93 & 281.76 & 272.87 & 266.37 & 256.82 \\
        \hline
    \end{tabular}
    \label{tab:Alaska data}
\end{table}

 \begin{figure}[H]
    \centering
    % Row for wind dataset
    {\includegraphics[height=6cm,width=0.45\textwidth]{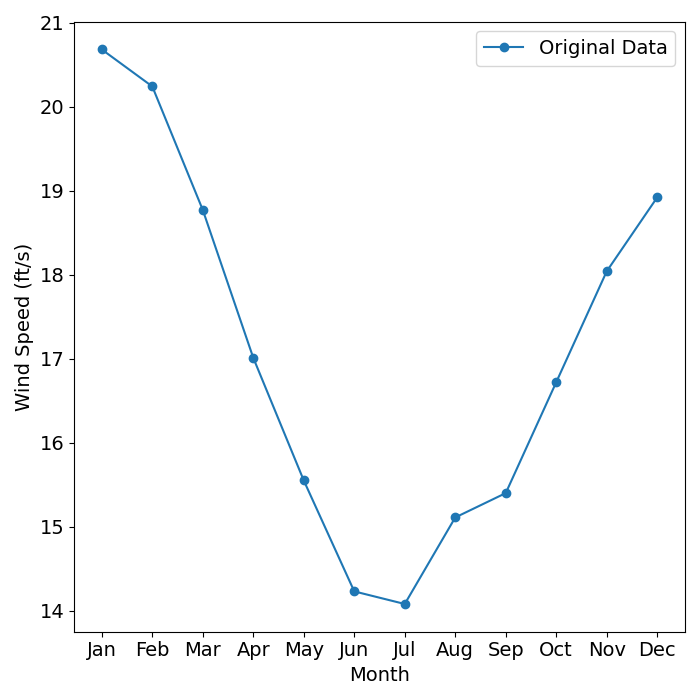}}
    {\includegraphics[height=6cm,width=0.45\textwidth]{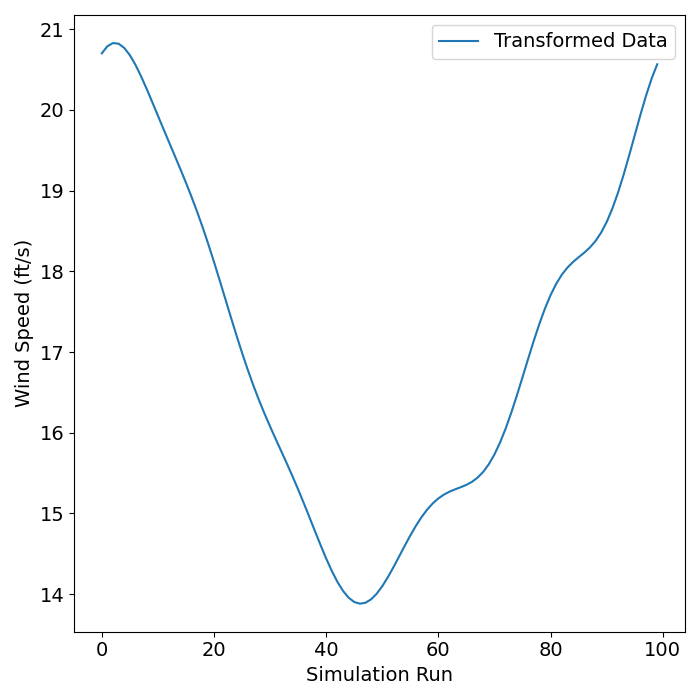}}
   \hfill
    {\includegraphics[height=6cm,width=0.45\textwidth]{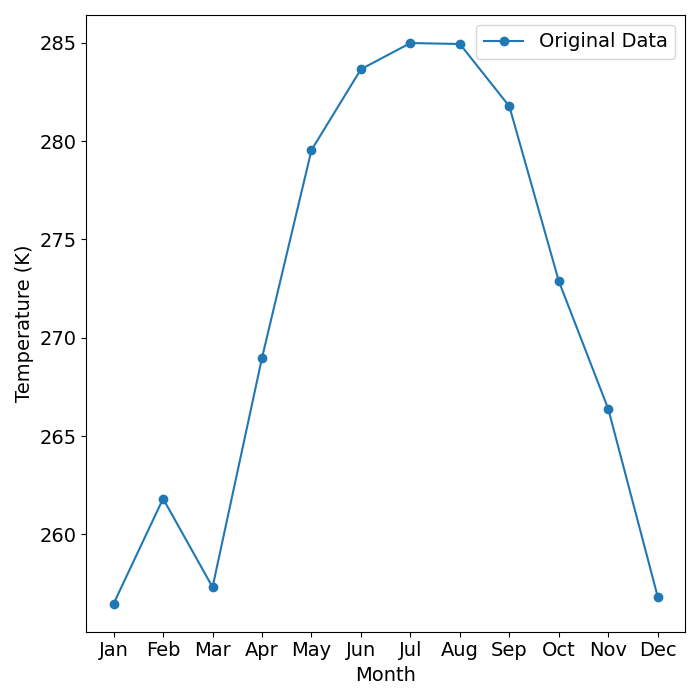}}
    {\includegraphics[height=6cm,width=0.45\textwidth]{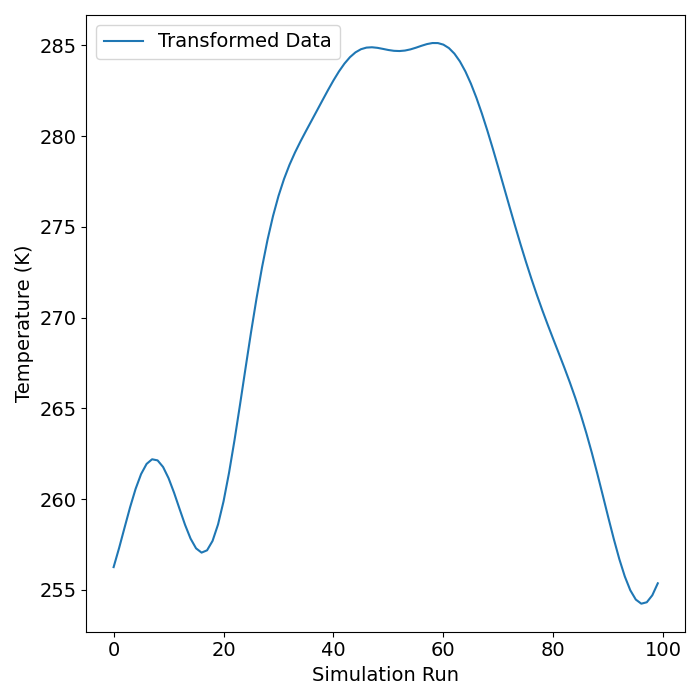}}
    \caption{Discrete and continuous variation of wind and temperature for Bethel, AK}
    \label{fig:Third}
   \end{figure} 
 
\end{itemize}

\subsection{Reliability of Transmission lines}

In this study, we focus on evaluating the reliability of transmission lines using a reliability model based on a limit state function, $g(R, S; t)$, where $R$ represents the critical thresholds, $S$ represents the maximum values of a dependent variable and $t$ represents time. Specifically, our analysis addresses reliability in terms of the maximum temperature or damage that a transmission line can withstand before failure. We define $S$ as the maximum temperature or damage of the transmission line, denoted by $\theta_{\text{max}}$ or $\varphi_{\text{max}}$, at any point in time $t$. We consider $R$ as the critical temperature or damage, represented as $\theta_{\text{lim}}$ or $\varphi_{\text{lim}}$, beyond which the integrity of the transmission line is compromised. Therefore, the limit state function is formulated as follows:
\begin{equation}
g(\theta_{\text{lim}}, \theta_{\text{max}}; t) = \theta_{\text{lim}} - \theta_{\text{max}}(t).
\end{equation}
\begin{equation}
g(\varphi_{\text{lim}}, \varphi_{\text{max}}; t) = \varphi_{\text{lim}} - \varphi_{\text{max}}(t).
\end{equation}

In general reliability analysis, $R$ and $S$ are usually considered as random variables. However, in our analysis, the thresholds $\theta_{\text{lim}}$ or $\varphi_{\text{lim}}$ are held as constants. The variability is fully attributed to $\theta_{\text{max}}(t)$ and $\varphi_{\text{max}}(t)$, which are functions of time and space. Consequently, the probability of failure $p_f(t)$ can be defined as

\begin{equation}
p_f(t) = P\{g(\theta_{\text{lim}}, \theta_{\text{max}}; t) < 0\}.
\end{equation}
\begin{equation}
p_f(t) = P\{g(\varphi_{\text{lim}}, \varphi_{\text{max}}; t) < 0\}.
\end{equation}

The thermo-electro-mechanical model introduced in this paper acts as a computational framework for assessing the operating temperature and damage evolution of the transmission line over its lifespan. Through the application of the Probabilistic Collocation Method (PCM), we subsequently estimate the probability of failure for the transmission line. 

\section{Methodology}
\subsection{Thermo-Electro-Mechanical Damage Phase-Field Model}
\begin{figure}[H]
    \centering
    \includegraphics[width=0.5\textwidth]{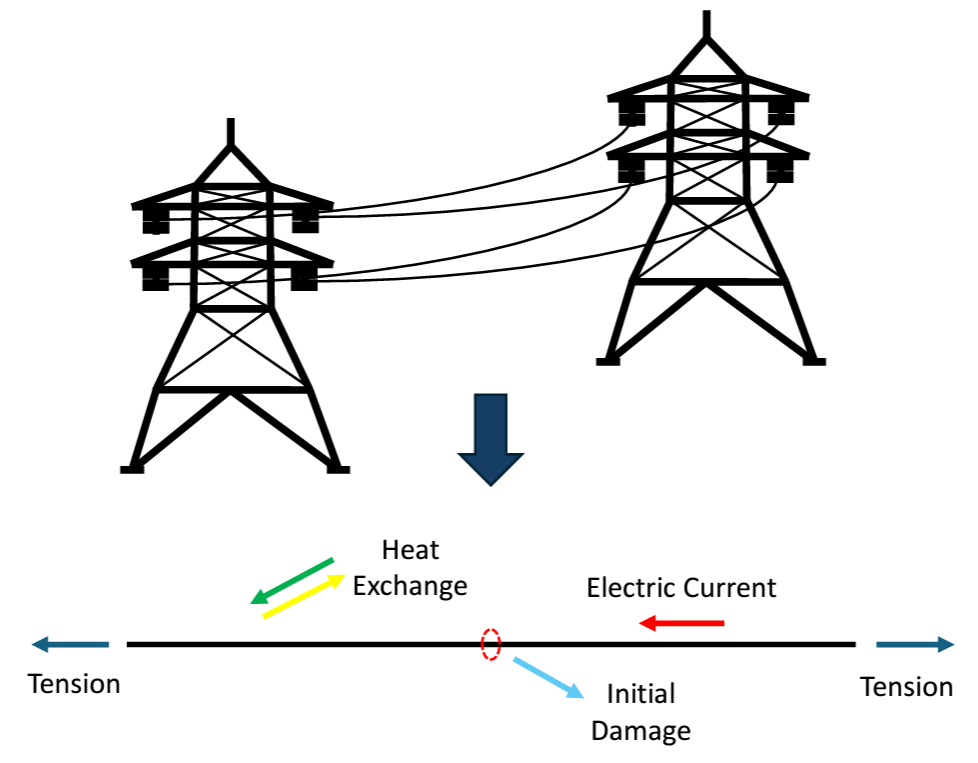}
    \caption{\textit{A one-dimensional representation of an initially damaged transmission line subjected to
mechanical tension, heat exchange with the environment, and an electric current.}}
    \label{fig:combine}
\end{figure}
In our study, we focus on a particular segment of a transmission cable. The segment is supported by two towers at the end and tensioned to ensure that sag remains within acceptable limits. We model the cable as a one-dimensional problem under tension. We consider its effective length as the projected span, as depicted in Fig.\ref{fig:combine}. While sag is an important factor in the dynamics of the cable, our analysis does not explicitly model it. Instead, we concentrate on the influence of the horizontal tension due to sag and assess their impact on material fatigue and damage over time.
The tension in the transmission cable is significantly affected by changes in its operating temperature. As the temperature rises, the cable undergoes thermal elongation, leading to an increased sag and a decrease in the horizontal tension. Conversely, a decrease in temperature causes the cable to contract, thus increasing the tension. The cyclic elongation and contraction due to sagging, over time, contributes to fatigue in the conductor leading to mechanical failure due to the accumulation of damage.
The thermal state of the conductor is determined by Joule heating, due to electrical currents, and is moderated by convective cooling due to the wind. Additionally, any damage in the conductor can lead to an increase in resistivity, resulting in further increase in the temperature. The increase in temperature accelerates the material's degradation and aging process. This dynamic highlights the interconnection between temperature and damage level in the conductor. The thermal model further incorporates heat transfer in the form of conduction, phase change, and radiation based on specific scenarios.  
The electrical model focuses mainly on the effects of the electric current flowing through a single cable on its material properties. The model does not consider the voltage levels supplied to consumers. As a result, resistive losses along the line increase the temperature of the cable, causing a voltage drop from one end to the other. For simplification, the overall electrical design of the transmission lines is not considered.

Overall, the proposed model suggests that any damage to the material increases its electrical resistance, which further increases the temperature-induced resistance and damage level. This leads to a positive feedback loop where damage accelerates the heating, and in turn, the increased temperature further accelerates the damage. This cycle can potentially lead to an early failure due to the thermal runaway and mechanical damage. To simulate these processes over a long time, we considered a quasi-static modeling approach that assumes the system reaches an equilibrium state more rapidly than the variation in the external loading conditions.

\subsubsection{Mechanical Model}
We use a non-isothermal phase-field framework to analyze damage and fatigue, based on the principles in \cite{boldrini2016non}. The framework consists of two partial differential equations (PDEs) describing the evolution of displacement ($u$) and damage ($\varphi$), accompanied by an ordinary differential equation (ODE) for the fatigue field ($\mathcal{F}$). The damage phase field $\varphi$ is representative of the volumetric fraction of the degraded material; with $\varphi = 0$ representing undamaged material and $\varphi = 1$ representing complete damage. Intermediate damaged states with $0 < \varphi < 1$ indicate varying degrees of material damage. The evolution of damage follows an Allen-Cahn-type equation, derived in conjunction with the equilibrium equation for the displacement field $u$, through the principle of virtual power and entropy inequalities, ensuring thermodynamic consistency. The fatigue field, denoted as $\mathcal{F}$, is considered an internal variable within the model, whose evolution is determined by constitutive relations that satisfy the entropy inequality across all permissible processes.

In our study, we employ a one-dimensional (1-D) representation for the mechanical body, such that it occupies the $\Omega \in \mathbb{R}$ at time $t \in (0,T]$. From the general governing equations, specific formulations of material evolution are derived based on the selection of free-energy potentials. Here, we consider the following free-energy function:
\begin{equation}
\label{eq:free_energy}
\Psi (\nabla u, \varphi, \nabla \varphi ,\mathcal{F}) 
= d(\varphi) Y \left(\nabla u\right)^2 + 
g_c \frac{\gamma}{2} \left(\nabla \varphi\right)^2  
+ \mathcal{K} (\varphi, \mathcal{F}) ,
\end{equation}
where $Y$ is the Young's modulus, $g_c$, is the fracture energy, and $\gamma > 0$ is the phase-field layer width parameter. The degradation function is considered as $d(\varphi) = (1-\varphi)^2$, not only affecting the elastic response, but also the electrical conductivity of the material beyond the elastic limit. Finally, $\mathcal{K} (\varphi, \mathcal{F})$ describes the evolution of the material damage in response to fatigue or aging.

The original model is considered as time-dependent in \cite{boldrini2016non}. However, we are interested in the long-term behavior of the material, assuming it reaches equilibrium between consecutive time steps allowing us to simplify the governing equations for $u$ and $\varphi$ into a quasi-static form. Yet, the evolution of $\mathcal{F}$ continues to be modeled as a time-dependent ODE as fatigue accumulates over time.

\begin{equation}
\label{eq:mechanical}
\nabla \cdot \,\left( (1-\varphi)^2 Y \nabla u \right) 
-  \gamma g_c \, \nabla \cdot \, ( \nabla \varphi \otimes \nabla \varphi  )+  f = 0,
\end{equation}

\begin{equation}
\label{eq:damage}
\gamma g_c \Delta \varphi  
+  (1- \varphi) (\nabla u)^T Y (\nabla u)
- \frac{1}{\gamma}  [ g_c \mathcal{H}' (\varphi) + \mathcal{F} \mathcal{H}_f' (\varphi) ] = 0,
\end{equation}

\begin{equation}
\label{eq:fatigue}
\dot{\mathcal{F}} = -  \frac{ \hat{F}}{\gamma}   \mathcal{H}_f (\varphi),
\end{equation}

The equations (\ref{eq:mechanical}) and (\ref{eq:damage}) are subjected to specific boundary conditions. 
The potentials $\mathcal{H}(\varphi)$ and $\mathcal{H}_f(\varphi)$ describe the evolution of damage from an undamaged state $(\varphi = 0)$ to a fully damaged state $(\varphi = 1)$ as the fatigue progresses from zero to a critical value $g_c$. By taking derivatives of these potentials with respect to $\varphi$, we obtain $\mathcal{H}'(\varphi)$ and $\mathcal{H}_f'(\varphi)$. To ensure a continuous and monotonically increasing transition, suitable choices for these potentials are:

\noindent\begin{minipage}{.55\linewidth}
	\begin{equation}
	\label{eq:potentials}
	\mathcal{H} (\varphi) =
	\begin{cases}
	0.5 \varphi^2   & \mbox{for} \; 0 \leq \varphi \leq 1 ,
	\vspace{0.1cm} 
	\\
	0.5 + \delta (\varphi -1)  &  \mbox{for} \; \varphi > 1 ,
	\vspace{0.1cm}
	\\
	- \delta \varphi &   \mbox{for} \; \varphi < 0 .
	\end{cases}
	\end{equation} 	
\end{minipage}%
\begin{minipage}{.45\linewidth}
	\begin{equation}
	\label{eq:potentials2}
	\mathcal{H}_f (\varphi) =
	\begin{cases}
	- \varphi & \mbox{for} \; 0 \leq \varphi \leq 1,
	\vspace{0.1cm} 
	\\
	-1 & \mbox{for} \; \varphi > 1,
	\vspace{0.1cm}
	\\
	\hspace{0.25cm} 0 & \mbox{for} \;  \varphi < 0.
	\end{cases}
	\end{equation}
\end{minipage}

The evolution of the fatigue field $\mathcal{F}$ is described through $\hat{F}$, which captures the formation and growth of micro-cracks in the material due to cyclic loading and temperature effects. The $\hat{F}$ is defined as a linear function of the stress levels associated with virgin materials:

\begin{equation} 
\label{eq:fhat}
\hat{F} =  \rho a \left(\frac{\theta_c}{\theta_0}\right) (1-\varphi) \left|Y \nabla u\right|,
\end{equation}
where the parameter $a$ is the aging rate, defined as the ratio of the conductor temperature $\theta_c$ to a reference temperature $\theta_0$ and $\rho$ is the material density.

The mechanical model considered here allows damage healing when tensile stress decreases. To avoid damage healing and irreversible damage process, we take an approach that considers a variable $\mathbb{H}$, representing the local maximum strain energy history, similar to that in \cite{miehe2010phase}. 

\begin{equation}
\label{eq:history}
\mathbb{H}(x,t) = \max((\nabla u(x,t))^T Y (\nabla u(x,t)),\mathcal{H}(x,t)).
\end{equation}

Incorporating the variable $\mathbb{H}$ into the damage equation, we obtain a new equation:
\begin{equation}
\label{eq:damage_H}
\gamma g_c \Delta \varphi  
+  (1- \varphi) \mathbb{H}
- \frac{1}{\gamma}  [ g_c \mathcal{H}' (\varphi) + \mathcal{F} \mathcal{H}_f' (\varphi) ] = 0.
\end{equation}

\subsubsection{Thermal Model}
The original model by \cite{boldrini2016non} relates fatigue to temperature increase due to repetitive and rapid loading cycles. However, our focus is on long-term material damage rather than short-term effects. To simplify the analysis and focus on the factors that significantly impact long-term performance, a quasi-static regime is assumed. This enables us to use the steady-state heat equation, to model the thermal behaviors under constant or slowly varying conditions. 
The steady-state heat equation is given by:
\begin{equation}
\label{eq:temperature}
\nabla \cdot (\kappa \nabla \theta) + q = 0.
\end{equation}
where $\kappa$ is the thermal conductivity, and 
$q$ is the net heat exchange. The heat exchange with the environment varies according to the different scenarios. To account for the effects of seasonal high winds and temperature, the following relation is adopted:

\begin{equation}
q = q_j - q_c,
\end{equation}
where $q_j$ is the heat due to the Joules heating and $q_c$ is a convective heat loss due to the wind. 

The convective heat loss due to cross flow over the cylinder \cite{cengel2011heat} is considered given by the relation: 

\begin{equation}
q_c = h (\theta_c - \theta_a),
\end{equation}
where $h$ is the convective heat transfer coefficient related to the Nusselt number $Nu_D$ and Prandtl number $Pr$ by the following relation: 
 \begin{equation}
Nu_D = \frac{hD}{\kappa_{air}} = C Re_D^m Pr^{\frac{1}{3}}
\end{equation}

We determine the Reynolds number $Re_D$ using a general relation: 
\begin{equation}
Re_D = \frac{vD}{\nu},
\end{equation}
where $\nu$ represents kinematic viscosity of air, $v$ represents the velocity of air, and $D$ represent the effective diameter of conductor. 

The experimental values of $C$ and $m$ are given in Table \ref{tab:experimental values} for different ranges of $Re_D$.\\
\begin{table}[h]
\caption{Values of \(C\) and \(m\) for different \(Re_D\) ranges}
\centering
\begin{tabular}{lll}
\hline
\(Re_D\) Range & \(C\) & \(m\) \\ \hline
0.4 -- 4 & 0.989 & 0.330 \\ 
4 -- 40 & 0.911 & 0.385 \\ 
40 -- 4000 & 0.683 & 0.466 \\ 
4000 -- 40,000 & 0.193 & 0.618 \\ 
40,000 -- 400,000 & 0.027 & 0.805 \\ \hline
\end{tabular}
\label{tab:experimental values}
\end{table}

Similarly, to account for the radiative effect of the wildfire, we consider the following relation given by \cite{guo2018determination} since the operating temperature of a transmission line is typically 5–15°C above the ambient temperature while the flame temperature of a wildfire can be up to 1200°C \cite{douglass2014real}.
\begin{equation}
q = q_j + q_f - q_c,
\end{equation}
where $q_f$ is the radiative heat transfer absorbed by the conductor:
\begin{equation}
\ q_f{} = \epsilon \sigma {\theta_f}^4 \phi \tau,
\end{equation}
where $\epsilon$ is flame emissivity, $\sigma$ is Stefan–Boltzmann constant, $\theta_f$ is the average temperature of the equivalent radiation surface, $\phi$ is the view factor, and $\tau$ is the atmospheric transmissivity.

We consider the effect of uniform layer of ice on the transmission line by the following relation:
\begin{equation}
q = q_j - q_i - q_c,
\end{equation}
where $q_i$ is the heat loss in terms of conduction and phase change given by the relation:
\begin{equation}
q_i = \frac{2 \pi  \kappa_{ice} \left(\theta_c - \theta_i\right)}{\ln \left(\frac{r_2}{r_1}\right)} + m_e  L_f,
\end{equation}
where $k_{ice}$ is the conductivity of ice, $r_1$ is the radius of cable, $r_2$ is the outer radius of cable with ice, $m_e$ is the mass of melted ice, and $L_f$ is the latent heat of fusion. There are several studies on the melting of ice over the transmission lines \cite{sadov2007mathematical,jiang2010simulation}, however, here we are only concerned about the heat loss due to the ice layer on the surface of the transmission line 

\subsubsection{Electrical Model}
The design of the transmission lines normally considers all the factors that influence the current flow $I_b$. However, here, we focus on understanding heating effects on the conductor due to Joule heating in the presence of initial damage, rather than on the overall design aspects of transmission lines. We therefore consider the current $I_b$ as an input to our model. Here we parameterize the allowable ampacity value as input current considering an all-Aluminum conductor of diameter 40 mm.

Additionally, we concentrate on the effect of Joules heating 
without involving the detailed time simulations of transient effects. To simplify our model efficiency, we assume the electric current as a DC-equivalent mean current, which remains constant over consecutive time steps allowing us to effectively measure quantities of interest.
This simplifies solving the conservation of current:
\begin{equation}
\nabla \cdot J = 0,
\end{equation}

\begin{equation}
J = \sigma_{E} E,
\end{equation}

\begin{equation}
E = -\nabla V,
\end{equation}
where $J$ is the electric current density, $E$ is the electric field, $V $ is the voltage, and  $\sigma_E$ is the electric conductivity of the wire which is the function of non-degraded conductivity through the degradation function $d(\varphi)$:
\begin{equation}
\sigma_E = (1 - \varphi)^2 \sigma_{E,T}.
\end{equation}

The non-degraded conductivity $\sigma_{E,T}$ is obtained from the reference temperature conductivity $\sigma_{E,0}$.
\begin{equation}
\sigma_{E,T} = \frac{\sigma_{E,0}}{1 + \alpha(\theta_c - \theta_0)},
\end{equation}
where $\alpha$ is the coefficient of resistivity of the conductor. To account for the effect of ice on the conductivity of the wire, we calculate the equivalent resistance by considering the wire and ice layer in parallel.

Finally, we obtain a partial differential equation for the voltage field using the above relations:
\begin{equation}
\nabla \cdot (-\sigma_{E} \nabla V) = 0,
\label{eq:voltage}
\end{equation}
with either $V$ or $J$ specified at the boundaries.

The Joule heating due to the current flowing through the wire is then defined as:

\begin{equation}
q_J = J \cdot E.
\end{equation}

In essence, as damage accumulates, it will lead to an increase in the voltage drop across the conductor due to the higher electrical resistance of the damaged conductor material. The increased resistance results in greater losses of electrical power through resistive heating effects, which further exacerbates the thermal loading on the conductor.

\subsubsection{Sag Consideration}

In this section, we consider, the effect of temperature variation on the horizontal tension acting on the cable supported by two towers, forming a catenary curve. The purely mechanical behaviors of the cables have been studied in  \cite{karoumi1999some,stengel2014finite}. However, we adopt a 1-D damage phase-field model, assuming the cable length $L$ is approximately equal to the span distance $S_L$ to simplify the analysis.  This assumption does not ignore the presence of sag $S$ but emphasizes that mechanical damage and fatigue are primarily driven by the horizontal tension component $H$.

We consider the effects of sag, especially under variable temperature conditions. High temperatures increase the length of the cable, increasing the sag and thereby reducing the horizontal tension. On the other hand, low temperatures 
 specifically in cold states, contract the cable, decreasing the sag and resulting in an increased tension. Thus, the mechanical load at one end of the cable is influenced by the initial pre-tension and the tension changes due to the operating temperature. The model determines the appropriate mechanical loading condition accordingly.

We follow the formulations outlined in \cite{grigsby2006electric} and consider $w_c$ as the weight per unit length and $H_0$ as the initial pre-tension, which is usually prescribed at about $20\%$ of the material's ultimate strength, to compute the initial sag $S_0$. The initial sag $S_0$ is defined as:

\begin{equation}
S_0 = \frac{w_c S_L^2}{8 H_0},
\end{equation}

For simplification, we assume the length of the cable $L$ is equal to the span length $S_L$, despite the presence of sag  $S$. Although simplified, the theoretical length required to accommodate the cable with sag over the span is given by:

\begin{equation}
L_0 = S_L + \frac{8 S_0^2}{3 S_L}.
\end{equation}

The change in cable length due to temperature variation is considered using the classical relation:

\begin{equation}
L = L_0 (1 + \alpha_L \Delta \theta),
\end{equation}
where $\alpha_L$ is the coefficient of thermal expansion. 

The resulting change in length is then used to determine the sag:

\begin{equation}
S = \sqrt{\frac{3 S_L (L - S_L)}{8}},
\end{equation}

Finally, the new sag is considered to determine the horizontal tension in the cable:

\begin{equation}
H = \frac{W S_L^2}{8S},
\end{equation}
where $W$ denotes the total weight, accounting for additional factors such as ice and wind along with the cable weight per unit length. 

\begin{equation}
W = \sqrt{(w_c+w_i)^2 + w_w^2},
\end{equation}

\begin{figure}[H]
    \centering
    \includegraphics[width=0.45\textwidth]{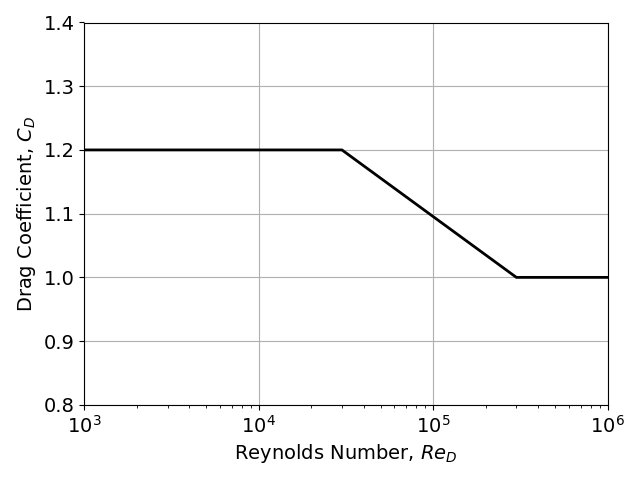}
    \caption{Drag coefficient $C_D$ for the transmission line.}
    \label{fig:CD}
\end{figure}

For the wind component, we consider the approach given by  \cite{reinoso2020wind,holmes2007wind} and calculate:
\begin{equation}
P_w = \frac{1}{2}\rho_{air} v^2,
\end{equation}
where $P_w$ is the wind pressure, $v$ is the wind velocity and $w_w$ is wind loading on the cable induced due to pressure:
\begin{equation}
w_w =  P_w C_D D \sin^2(\theta_w) \alpha,
\end{equation}
where $\rho_{air}$ represents the density of air, $D$ is the diameter, $\theta_w$ represents the angle between the line and wind flow, $\alpha$ is the span factor, which we consider urban terrain and $C_D$ is the drag Coefficient. The drag coefficient is obtained using the relation with Reynold's number as shown in Fig.\ref{fig:CD}.

To account for the ice load on the transmission cable, we follow:
\begin{equation}
    V = \pi (D + t_{ice}) t_{ice},
\end{equation}
\begin{equation}
    w_i = \rho_{i} V g, 
\end{equation}
where $D$ is the diameter of the conductor, $t_{ice}$ is the ice thickness, $\rho_i$ is the density of ice, $V$ is the volume of ice per unit length, and $g$ is the acceleration due to gravity.

The wind loading under icing is considered using the following relation:
\begin{equation}
w_w =  P_w C_D (D + 2t) \sin^2(\theta) \alpha
\end{equation}

\subsubsection{Multi-physics framework}

\begin{figure}[H]
    \centering
    \includegraphics[width=0.5\textwidth]{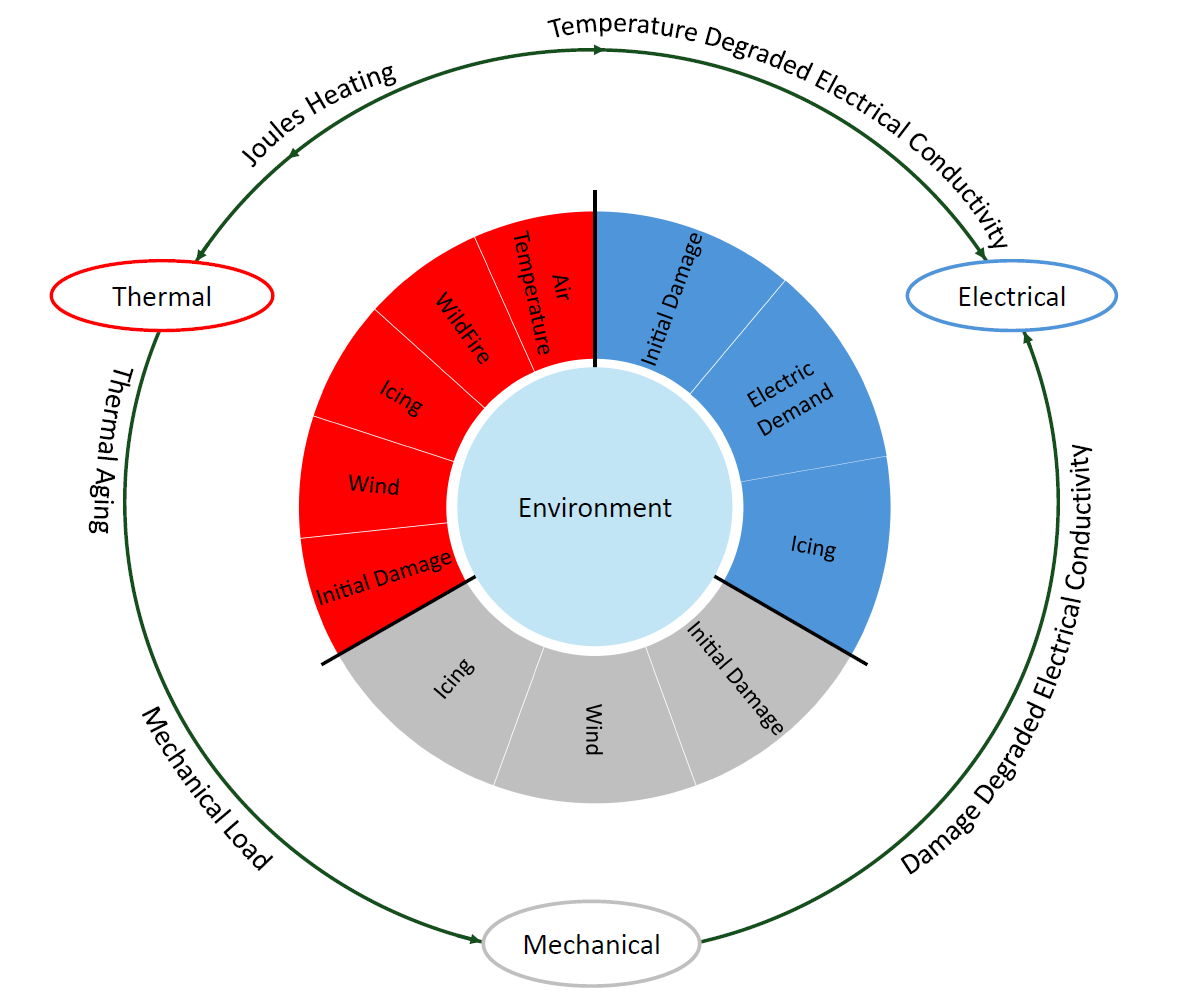}
    \caption{\textit{Schematic diagram representing the relationship between the thermal, electrical, and mechanical models, in addition to an abstract, scenario-dependent environmental module that provides input and initial conditions.}}
    \label{fig:relation}
\end{figure}

We illustrate a schematic diagram for modeling the transmission line failure in Fig~\ref{fig:relation}. The main governing equations for the coupled system are given by Eqs. (\ref{eq:mechanical}), (\ref{eq:damage_H}), (\ref{eq:fatigue}), (\ref{eq:temperature}), and (\ref{eq:voltage}), along with equations for fatigue potentials, thermal and mechanical loads, and the degradation of electrical conductivity. The diagram provides a comprehensive view of the interactions between the thermal, electrical, and mechanical components. Additionally, we introduce an abstract representation of an Environmental module, which determines the inputs and initial conditions for all other modules based on the specific scenario being considered.

\section{Deterministic Solution}
\label{sec:deterministic}

In this section, we discuss the Finite-Element discretization of the proposed multi-physics framework and describe the solution procedure for the deterministic case. The deterministic solution acts as a black-box for non-intrusive stochastic methods used in uncertainty quantification and serves as a clear guide for further evaluations of multi-dimensional uncertainty propagation.

\subsection{Finite-Element Discretization}
We consider a one-dimensional domain of length $L$ = 200 $m$ for solving the problem using the finite element method. We used the linear shape function to obtain the spatial discretization considering the number of elements $N$ = 1000. The governing equations, initially in volumetric form, are multiplied by the cross-sectional area $A$ to obtain a one-dimensional representation. Then, the equations are multiplied with a test function $w$, to obtain the weak forms after performing integration by parts. 

\begin{equation}
\int_0^L -(1 - \varphi)^2 Y A(x) \frac{du}{dx}\frac{dw}{dx} dx + \int_0^L \gamma g_c A(x) \left(\frac{d\varphi}{dx}\right)^2 \frac{dw}{dx} dx + \int_0^L f A(x) w dx = 0,
\end{equation}

\begin{equation}
\begin{split}
\int_0^L - \gamma g_c A(x) \frac{d\varphi}{dx}\frac{dw}{dx} dx + \int_0^L A(x) \mathbb{H} w dx - \int_0^L A \mathbb{H} \varphi w dx \\ - \int_0^L \frac{g_c A(x)}{\gamma} \varphi w dx + \int_0^L \frac{A(x)}{\gamma}\mathcal{F}w dx = 0,
\end{split}
\end{equation}

\begin{equation}
\int_0^L \dot{\mathcal{F}} w A(x) dx = \int_0^L \frac{-\rho a (1 - \varphi) Y \vert\frac{du}{dx}\vert (-\varphi)}{\gamma}\frac{\theta_c}{\theta_0} w A(x) dx, 
\end{equation}

\begin{equation}
\begin{split}
\int_0^L - \kappa A(x) \frac{d\theta}{dx}\frac{dw}{dx} dx + \int_0^L \sigma_E A(x)\left(\frac{dV}{dx}\right)^2 w dx \\- \int_0^L h \theta A_s(x) w dx + \int_0^L h \theta_{a} A_s(x) w dx = 0,
\end{split}
\end{equation}

\begin{equation}
\int_0^L \sigma_E \frac{dV}{dx} \frac{dw}{dx} A(x) w dx = 0,
\end{equation}
where we define $A_s(x)$ as the variable surface area where convective heat transfer due to wind occurs.

In the case of wildfire, the weak form of the heat equation is represented as:
\begin{equation}
\begin{split}
\int_0^L - \kappa A(x) \frac{d\theta}{dx}\frac{dw}{dx} dx + \int_0^L \sigma_E A(x)\left(\frac{dV}{dx}\right)^2 w dx \\- \int_0^L h \theta A_s(x) w dx + \int_0^L h \theta_{a} A_s(x) w dx + \int_0^L \epsilon \sigma {\theta_f}^4 \phi \tau A_s(x) w dx = 0,
\end{split}
\end{equation}

Similarly, in the case of icing, we obtain the weak form as follows:
\begin{equation}
\begin{split}
\int_0^L - \kappa A(x) \frac{d\theta}{dx}\frac{dw}{dx} dx + \int_0^L \sigma_E A(x)\left(\frac{dV}{dx}\right)^2 w dx \\- \int_0^L \frac{2 \pi  k \theta_c} {\ln(\frac{r_2}{r_1})} w dx + \int_0^L \frac{2 \pi  k \theta_i} {\ln(\frac{r_2}{r_1})} w dx + \int_0^L \rho_{i} \pi (D + t) t  L_f w dx = 0,
\end{split}
\end{equation}
where convective heat transfer is not considered, as the ice layer on the transmission line prevents the direct contact with the wind.

We adopt a linear approximation for each element $k$, where the field variables are expressed as a linear combination of nodal basis functions:

\begin{align}
u^k &= N \hat{u}^k,\\
\varphi^k &= N \hat{\varphi}^k ,\\
\mathcal{F}^k &= N \hat{\mathcal{F}}^k ,\\ 
\theta^k &= N \hat{\theta}^k ,\\
V^k &= N \hat{V}^k.
\end{align}

Finite-element interpolations for spatial derivatives are computed as linear combinations of the derivatives of the shape functions:

\begin{align}
\left(\frac{du}{dx}\right)^k &= B \hat{u}^k,\\
\left(\frac{d\varphi}{dx}\right)^k &= B \hat{\varphi}^k,\\
\left(\frac{d\theta}{dx}\right)^k &= B \hat{\theta}^k,\\
\left(\frac{dV}{dx}\right)^k &= B \hat{V}^k,\\
\end{align}

\noindent where we define $N$, $B$, $\hat{u}^k$, $\hat{\varphi}^k$, $\hat{\mathcal{F}}^k$, $\hat{\theta}^k$, $\hat{V}^k$ as
\begin{align}
N &= \begin{bmatrix}
N_1 & N_2
\end{bmatrix},\\
B &= \begin{bmatrix}
N_{1,x} & N_{2,x}
\end{bmatrix},\\
\hat{u}^k &= \begin{bmatrix}
u_1^k & u_2^k
\end{bmatrix},\\
\hat{\varphi}^k &= \begin{bmatrix}
\varphi_1^k & \varphi_2^k
\end{bmatrix},\\
\hat{\mathcal{F}}^k &= \begin{bmatrix}
\mathcal{F}_1^k & \mathcal{F}_2^k
\end{bmatrix},\\
\hat{\theta}^k &= \begin{bmatrix}
\theta_1^k & \theta_2^k
\end{bmatrix},\\
\hat{V}^k &= \begin{bmatrix}
V_1^k & V_2^k,
\end{bmatrix}.
\end{align}
where $N_1$ and $N_2$ are linear interpolation functions.

We substitute the previous approximations into the weak form and use a forward Euler method for the evolution of $\mathcal{F}$, resulting in the following discretization for each $k$th element: 

\begin{align}
K_u \hat{u}^k &= w_u + M \hat{f}^k,\\
K_{\varphi} \hat{\varphi}^k &= w_\varphi,\\
M \hat{\mathcal{F}^{n+1}}^k &= M \hat{\mathcal{F}^{n}}^k + \Delta t w_{\mathcal{F}},\\
K_\theta \hat{\theta}^k &= w_\theta,\\
K_V \hat{V}^k &= 0.
\end{align}
where the superscripts $n$ and $n+1$ denote the current and next time steps, respectively. The discrete forms are defined using these operator conventions:
\begin{align}
K_u &= \int_k (1 - N\hat{\varphi}^k)^2 Y A(x) B^T B \, dx,\\
w_u &= \int_k \gamma g_c A(x) (B \hat{\varphi}^k)^2 B \, dx,\\
M &= \int_k A(x) N^T N \, dx,\\
K_\varphi &= \int_k \gamma g_c A(x) B^T B \, dx + \int_k \mathbb{H} A(x) N^T N \, dx + \int_k \frac{g_c A(x)}{\gamma} N^T N \, dx,\\
w_\varphi &= \int_k \mathbb{H} A(x) N \, dx + \int_k \frac{A(x)}{\gamma} N^T \hat{\mathcal{F}^n}^k N \, dx,\\
K_\theta &= \int_k \kappa A(x) B^T B \, dx + \int_k h A_s(x) N^T N \, dx,\\
w_\theta &= \int_k \sigma_E A(x) \left(B \hat{V}^k\right)^2 N \, dx  + \int_k h A_s(x) \theta_{a} N \, dx,\\
K_v &= \int_k (1 - N\hat{\varphi}^k)^2 \sigma_{E,T} A(x) B^T B \, dx.
\end{align}
To account the wildfire, $w_{\theta}$ is modified as:
\begin{align}
w_\theta &= \int_k \sigma_E A(x) \left(B \hat{V}^k\right)^2 N \, dx  + \int_k h A_s(x) \theta_{a} N \, dx + \int_k \epsilon \sigma {\theta_f}^4 \phi \tau A_s(x) N dx
\end{align}
while, in the case of icing, $K_{\theta}$ and $w_{\theta}$ are modified as:
\begin{align}
K_\theta &= \int_k \kappa A(x) B^T B \, dx + \int_k \frac{2 \pi  k \theta_c} {\ln(\frac{r_2}{r_1})} N^T N \, dx, \\
w_\theta &= \int_k \sigma_E A(x) \left(B \hat{V}^k\right)^2 N \, dx  + \int_k \frac{2 \pi  k \theta_i} {\ln(\frac{r_2}{r_1})} N dx - \int_0^L \rho_{i} \pi (D + t) t  L_f N dx
\end{align}

We then assemble the local matrices and vectors to obtain the respective global forms using standard finite-element assembly procedures and implement the two-point Gauss quadrature rule for the integration. 
\begin{align}
K_u \hat{u} &= w_u + M \hat{f},\\
K_{\varphi} \hat{\varphi} &= w_\varphi, \\
M \hat{\mathcal{F}}^{n+1} &= M \hat{\mathcal{F}}^n  + \Delta t w_{\mathcal{F}}, \\
K_\theta \hat{\theta}  &= w_\theta, \\
K_V \hat{V} &= 0.
\end{align}

The general algorithm for the deterministic solution is summarized in Algorithm~\ref{algo:deterministic}:
\begin{algorithm}[H]
\caption{Solution of Thermo-Electro-Mechanical Model.}
\label{algo:deterministic}
\begin{algorithmic}[1]
\State Choose initial pre-tension.
\For{Each time-step}
\State Compute the current tensile load.
\State Solve for displacements.
\State Update strain energy history.
\State Solve damage field.
\State Update fatigue.
\State Solve the temperature field.
\State Solve voltage field.
\EndFor
\end{algorithmic}
\end{algorithm}

\subsection{Numerical Results}

The reliability of transmission lines also depends upon the type of materials. The most common are Aluminum conductor steel-reinforced (ACSR), all-Aluminum alloy conductor (AAAC), and all-Aluminum conductor (AAC). Here, we consider the material properties of aluminum assuming all-Aluminum 
conductor-type transmission line. The conductor is subjected to cyclic loading related to the ambient temperature, wind, and current. The cyclic loading of wind and temperature for each specific scenario is detailed in \ref{sec:Problem Statement}. However, due to difficulty in assessing current loading data, we parameterize it using the following relation:
\begin{equation}
    I(t) = -I_b - I_a (\sin 4 \pi t),
\end{equation}
where $I_b$ is the base current and $I_a$ is the amplitude. The base value for the current is considered as 1500 A, representing the allowable ampacity for a 40 mm Aluminum conductor.

For the boundary conditions, we set the horizontal tension by specifying $u=0$ at $x=0$ and $H$ at $x=L$. Similarly, for damage, we use the Neumann boundary condition $\frac{d \varphi}{dx} = 0 $ at both boundaries. We enforce the Dirichlet boundary condition for the current conservation equation, similar to those in the mechanical case: setting $V = 0$ at $x=0$ and imposing a current density $J$ at $x=L$.

\begin{figure}[H]
    \centering
    \includegraphics[width=0.45\textwidth]{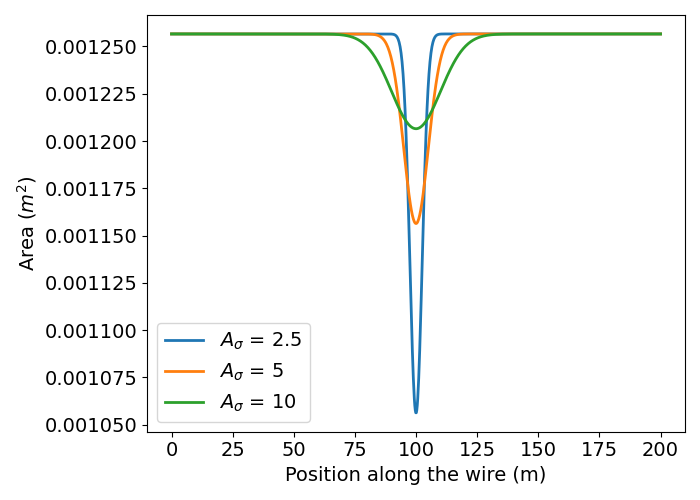}
    \caption{Cross-sectional area variation as a function of $A_\sigma$}
    \label{fig:Area1}
\end{figure}

In general, materials inherently possess imperfections that accumulate over time, eventually reaching a critical point where failure occurs. To model this, we consider the initial damage to be a variable cross-sectional area at the midpoint of the conductor, representing the overall effects of multiple defects that lead to fracture. For comparative purposes, we establish three scenarios where materials are assumed to have different sizes of initial damage. We define the cross-sectional area using the following relation:
\begin{equation}
A(x) = A_0 \left(1 - \frac{1}{A_\sigma \sqrt{2 \pi}} \exp\left( \frac{- (x - L/2)^2 }{2 A_\sigma^2}\right)\right),
\end{equation} 
where $A_0 $ represents the cross-sectional area of the conductor in its undamaged state, while $A_\sigma$ quantifies the ratio of spread to depth of the variation in the area, reflecting various degrees of damage. Fig.\ref{fig:Area1} provides illustrations of different area profiles corresponding to values $A_\sigma$.
The properties of aluminum conductor, air, and ice are provided in Table.\ref{tab:param}. All the simulations are conducted over a period of 5000 runs with time steps of $\delta t$ = 0.01, equivalent to a lifecycle of 50 years. Considering the behavior of Aluminum conductors, which begin to anneal at temperatures exceeding 366 K \cite{vasquez2017end,hathout2018impact}, and rupture at temperatures above 373 K \cite{cimini2013temperature}, we set 373 K as the maximum allowable temperature $\theta_{\text{lim}}$. We consider the maximum damage limit to $\varphi_{\text{lim}} = 0.8$. The simulations are stopped when the threshold for either of these parameters is reached.

\begin{table}[H]
\caption{Properties of Aluminum, Air, and Ice.}
    \centering
    \begin{tabular}{lll} % Simple table with vertical lines
        \hline
        \textbf{Property} & \textbf{Value} & \textbf{Unit} \\ 
        \hline
        Young modulus $Y$ & 69 & GPa \\
        Damage layer width $\gamma$ & 0.02 & m \\
        Fracture energy $g_c$ & 10 & kN/m \\
        Density $\rho$ & 2700 & kg/m$^3$ \\
        Aging coefficient $a$ & $1 \times 10^{-10}$ & m$^5$/(y kg) \\
        Thermal conductivity $\kappa$ & 237 & W/(m K) \\
        Electrical conductivity $\sigma_{E,0}$ & $3.77 \times 10^7$ & S/m \\
        Temperature coefficient $\alpha$ & $3.9 \times 10^{-3}$ & K$^{-1}$ \\
        \hline
        Density of air $\rho_{air}$ & 1.225 & kg/m$^3$ \\
        Kinematic viscosity of air $\nu$ & $15 \times 10^{-6}$  &  m$^2$/s\\
        Thermal Conductivity of air $\kappa_{air}$ & 0.0295  &  W/(m K)\\
        Prandtl Pr & 0.71 &  \\
        \hline
        Density of ice $\rho_{ice}$ & 917 & kg/m$^3$ \\
        Resistivity of ice & $1 \times 10^{9}$  &  $\Omega$ m\\
        Thermal Conductivity of ice $\kappa_{ice}$ & 2.39  &  W/(m K)\\
        Latent Heat of Fusion Le & $3.36 \times 10^{5}$ & J/kg  \\
        \hline
    \end{tabular}
    \label{tab:param}
\end{table}

\begin{figure}[H]
	\centering
	\subfloat[Damage.]{\includegraphics[width=0.25\textwidth]{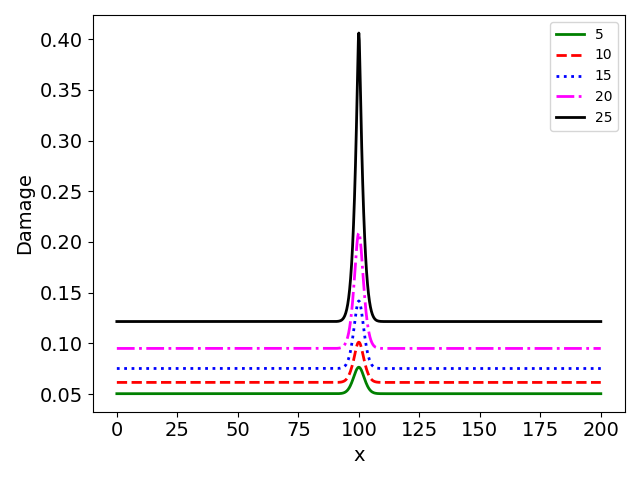}}
	\subfloat[Fatigue.]{\includegraphics[width=0.25\textwidth]{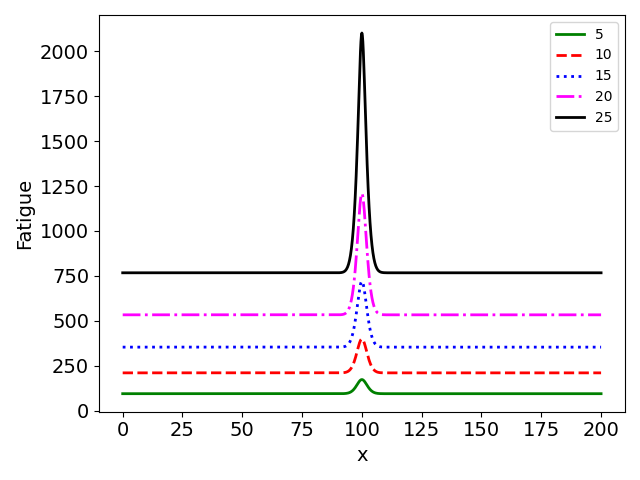}}
	\subfloat[Temperature.]{\includegraphics[width=0.25\textwidth]{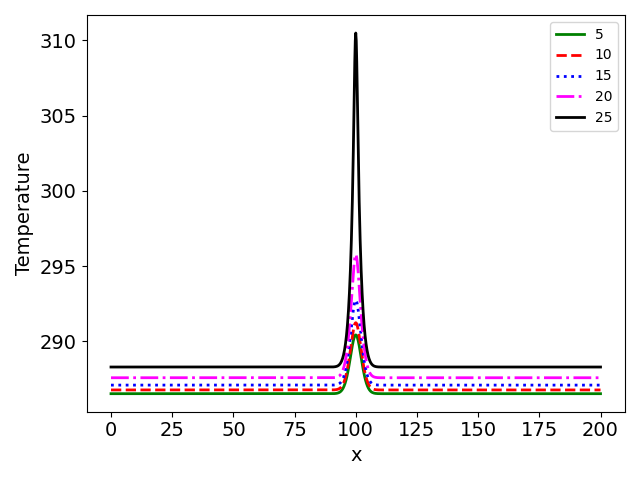}}
	\subfloat[Voltage drop.]{\includegraphics[width=0.25\textwidth]{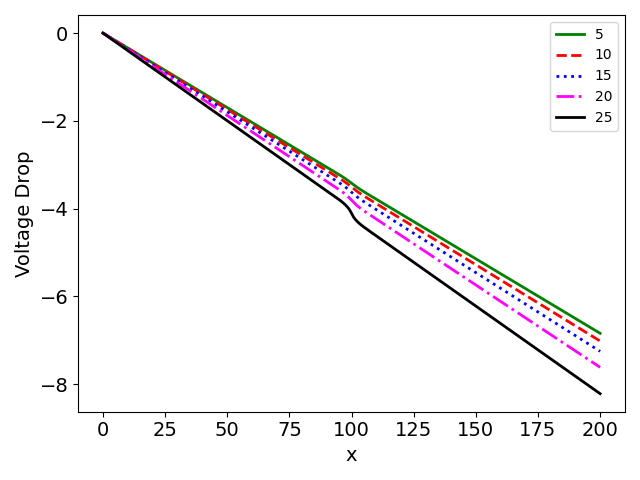}}
	\caption{Evolution of field variables over 5 years of interval.}
	\label{fig:field_1_1}
\end{figure}

We first analyze the evolution of field quantities using the properties presented in Table.\ref{tab:param}. We plot the progression in damage, fatigue, temperature, and voltage fields every 5 years as shown in Fig.\ref{fig:field_1_1}. We observe that damage typically originates and accumulates in regions with reduced cross-sectional areas, leading to elevated temperatures and noticeable distortions in the voltage fields. Furthermore, as damage progresses and temperatures rise, there is a corresponding increase in voltage drop along the line, due to the increased electrical resistance of the conductor.

\begin{figure}[H]
	\centering
	\subfloat[Initial damage.]{\includegraphics[width=0.45\textwidth]{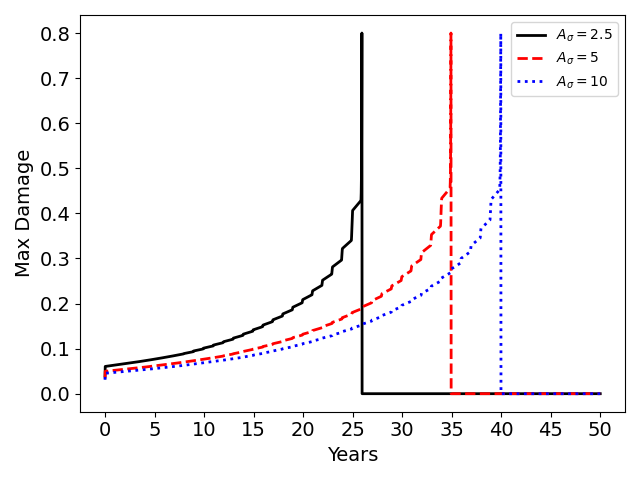}}
	\subfloat[Extreme wind.]{\includegraphics[width=0.45\textwidth]{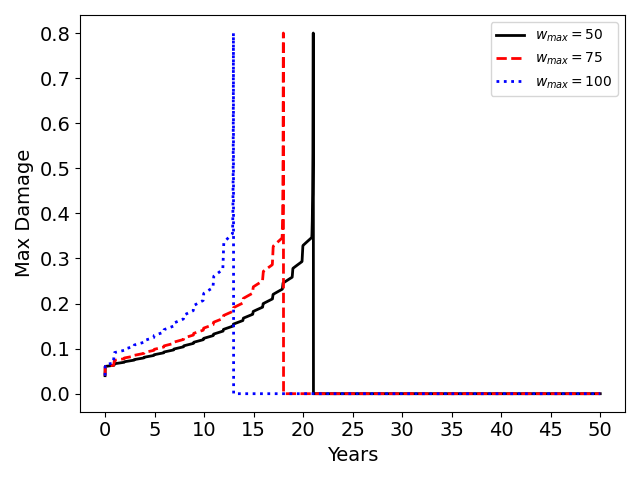}}
	\caption{Effects of initial damage and extreme wind parameter $w_{max}$ in Scenario 1.}
	\label{fig:Scenario 1}
\end{figure}

We continue our investigation by analyzing the failure of the transmission line using the properties and loading conditions. Initially, we vary the initial damage $A_{\sigma}$, to observe its impact on the transmission line's longevity. Then, we introduce the effect of high wind, wildfire, and icing to the respective scenarios. In Scenario 1, which is characterized by high wind conditions and variable temperatures across different seasons, we assess how these environmental factors exacerbate the existing damage leading to failure. By solely adjusting $A_{\sigma}$ as shown in Fig.\ref{fig:Scenario 1}, we notice a significant decrease in the lifespan of the line with more severe initial damage, highlighting the critical role of initial damage in failure processes. 

\begin{figure}[H]
	\centering
	\subfloat[Initial damage.]{\includegraphics[width=0.45\textwidth]{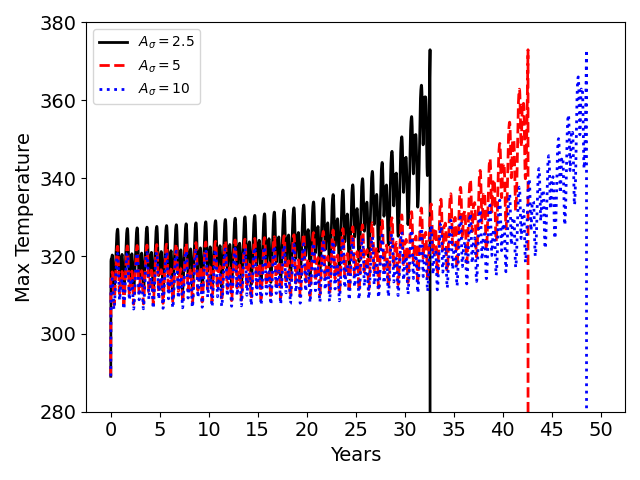}}
	\subfloat[Wildfire.]{\includegraphics[width=0.45\textwidth]{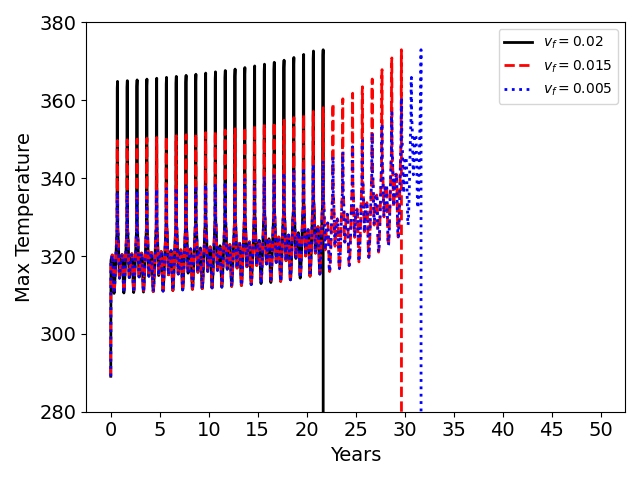}}
	\caption{Effects of initial damage and wildfire in Scenario 2.}
	\label{fig:Scenario 2}
\end{figure}

Subsequently, we introduce high seasonal wind speeds into the analysis of an already damaged transmission line to further explore the impact of unexpected high wind on the lifespan of the infrastructure. While high wind speeds facilitate convective cooling, the mechanical load imposed on the line ultimately exacerbates the damage over time. This pushes the damage level beyond the critical thresholds, failing the transmission line. Fig.\ref{fig:Scenario 1} shows the life span of lines reduces with an increase in wind speed. In Scenario 2, we adopt a similar procedure by varying the initial damage. However, the observed mode of failure here is the operating temperature, which is influenced by the region's high ambient temperatures and low wind speeds that reduce the convective cooling effect. Furthermore, to investigate the impact of wildfire in the form of radiative heat, we adjust the distance between the wildfire and the transmission line in the form of view factor. Fig.\ref{fig:Scenario 2} shows the shorter distance intensifies the effects of the radiative heat on the transmission line, leading to an early failure.

\begin{figure}[H]
	\centering
	\subfloat[Initial damage.]{\includegraphics[width=0.45\textwidth]{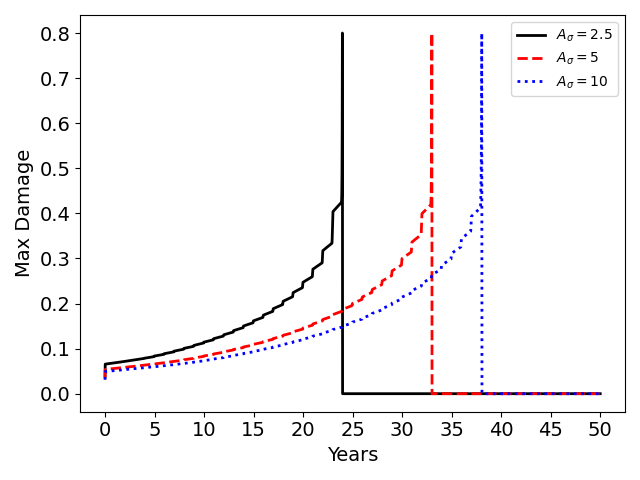}}
	\subfloat[Wildfire.]{\includegraphics[width=0.45\textwidth]{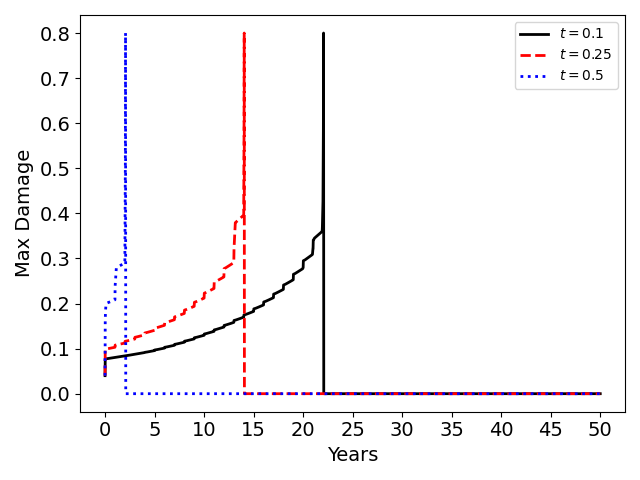}}
	\caption{Effects of initial damage and icing in Scenario 3.}
	\label{fig:Scenario 3}
\end{figure}

In Scenario 3, we explore the impact of icing on the transmission line in a region characterized by consistently low temperatures and high wind speed throughout the year. Initially, we assess the influence of initial damage and find that more severe initial damage leads to earlier failure due to cumulative damage effects as shown in Fig.\ref{fig:Scenario 3}. Subsequently, we adjust the thickness of the ice to examine the impact of ice loading on failure. We observe that increasing the ice thickness substantially shortens the lifespan of the transmission line due to the additional mechanical loading imposed on the line.

\section{Stochastic Solution}
For the stochastic solution, we treat the model as a black box for performing uncertainty quantification (UQ), sensitivity analysis (SA), and probability of failure ($p_f$). Non-intrusive methods such as PCM are particularly advantageous as they allow us to apply the same procedures from the deterministic solution to each realization of the stochastic analysis. Therefore, we utilize the PCM to compute the moments of our quantity of interest (QoI), which in this study is the conductor temperature or damage based on the specific scenario.

Additionally, PCM is used for global sensitivity analysis by computing the Sobol Sensitivity Index $S_i$. This index evaluates the contribution of input parameters to the total variance of the QoI solution. Using the PCM, the sensitivity index can be obtained just by post-processing UQ data in a straightforward fashion.

Lastly, the PCM framework is used to compute the probability of failure directly from the first moment of a Bernoulli random variable defined by $g(\theta_{lim},\theta_{max};t)$ or $g(\varphi_{lim},\varphi_{max};t)$

\subsection{Uncertainity Quantification}

To perform the UQ analysis, we employ PCM which uses polynomial interpolation to approximate solutions in the stochastic space. The PCM maps points from physical to stochastic space using parametric probability density functions (PDF). This is achieved by using orthogonal Lagrange polynomials to approximate the solution. Due to the orthogonality properties, the computation of expectation and variance to assess the quantity of interest (QoI) at the collocation points becomes straightforward. Using this approach, the computational cost is significantly reduced and the convergence rates are improved. 

Following the methodology presented in \cite{barros2021integrated}, we consider a complete probability space $(\Omega_s, \mathcal{G}, \mathbb{P})$, where $\Omega_s$ is the space of outcomes $\omega$, $\mathcal{G}$ is the $\sigma$-algebra, and $\mathbb{P}$ is a probability measure, $\mathbb{P}: \mathcal{G} \to [0,1]$. The transmission line model becomes stochastic by treating material and load parameters as random variables in a set $\xi(\omega)$, resulting in outputs as temperature or damage also being random variables. For simplicity, we represent the random parameters as $\xi = \xi(\omega)$.

We denote our quantity of interest by $Q$ and express its mathematical expectation $\mathbb{E}\left[ Q(x, t; \xi) \right]$ within a one-dimensional stochastic space as:
\begin{equation}
\label{eq:pcm_expectation}
\mathbb{E}\left[ Q(x, t; \xi) \right] = \int_{a}^b Q(x, t; \xi) \rho(\xi) d \xi,
\end{equation}
where $\rho(\xi)$ is the PDF of $\xi$. We perform the integration using Gauss quadrature, which maps the physical parametric domain to the standard domain $[-1, 1]$. The integral can then be expressed as:
\begin{equation}
\label{eq:pcm_expectation2}
\mathbb{E}\left[ Q(x, t; \xi) \right] = \int_{-1}^1 Q(x, t; \xi(\eta)) \rho(\xi(\eta)) J d \xi(\eta),
\end{equation}
where $J = d\xi / d\eta$ is the Jacobian of the transformation. We approximate the expectation using polynomial interpolation of the solution in the stochastic space, represented as  $\hat{Q}(x, t; \xi)$:
\begin{equation}
\label{eq:pcm_expectation3}
\mathbb{E}\left[ Q(x, t; \xi) \right] \approx \int_{-1}^1 \hat{Q}(x, y,t; \xi(\eta)) \rho(\xi(\eta)) J d \xi(\eta).
\end{equation} 
The solution is interpolated in the stochastic space using Lagrange polynomials $L_i(\xi)$:
\begin{equation}
\label{eq:poly_approximation}
\hat{Q}(x, t; \xi) = \sum_{i = 1}^{I} Q(x, t; \xi_i) L_i(\xi),
\end{equation}
which satisfies the Kronecker delta property at the interpolation points:
\begin{equation}
\label{eq:delta}
L_i(\xi_j) = \delta_{ij}.
\end{equation}

Substituting the polynomial approximation from Eq.(\ref{eq:poly_approximation}) into Eq.(\ref{eq:pcm_expectation3}), and using the quadrature rule, we approximate the integral and compute the expectation as

\begin{equation}
\mathbb{E}\left[ Q(x, t; \xi) \right] \approx \sum_{p = 1}^{P} w_p \rho(\xi(\eta)) J \sum_{i = 1}^{I} Q(x, t; \xi(\eta) ) L_i(\xi(\eta) ), \label{eq:pcm_exp_approx2}
\end{equation}
where we calculate the coordinates $\eta_p$ and weights $w_p$ for each integration point $q = 1,\, 2,\, \dots,\, P$. This is obtained efficiently by choosing the collocation points same as the integration points in the parametric space. Using the Kronecker property of the Lagrange polynomials, as detailed in Eq.(\ref{eq:delta}), we simplify the approximation from Eq.(\ref{eq:pcm_exp_approx2}) into a single summation:
\begin{equation}
\label{eq:pcm_single}
\mathbb{E}\left[ Q(x, t; \xi) \right] = \sum_{p = 1}^{P} w_p \rho(\xi_p(\eta_p))  J Q(x, t; \xi_p(\eta_p)).
\end{equation}

We use a linear affine mapping from the standard domain to the real domain using the formula $\xi_p(\eta_p) = a + \frac{(b-a)}{2}(\eta_p + 1)$. This yields the Jacobian for a one-dimensional integration as $J = (b-a)/2$ and determines the values of the random variable in the physical space.
Finally, we express the integration as a summation over the collocation points. Assuming a uniform distribution for the parameters in the physical space, where $\xi \sim \mathcal{U}[a,b]$ with $\rho(\xi) = 1/(b-a)$, the expectation becomes:
\begin{equation}
\mathbb{E}\left[ Q(x, y, t; \xi) \right] = \frac{1}{2}\sum_{p = 1}^{P}  w_p Q(x, t; \xi_p).
\end{equation} 

The standard deviation is computed as
\begin{equation}
\sigma \left[ Q(x, t; \xi) \right] =  \sqrt{ \frac{1}{2}\sum_{p = 1}^{P}  w_p \left( Q(x, t; \xi_p) - \mathbb{E}\left[ Q(x, t; \xi) \right] \right)^2}.
\end{equation}

Generalization of PCM to higher dimensions involves adding additional integrals to Eq.~(\ref{eq:pcm_expectation}) which then reduces to

\begin{align}
\mathbb{E}\left[ Q(x, t; \xi^1,\, \dots ,\, \xi^k) \right] &= \mathbb{E}_{PCM}\left[ Q(x, t; \xi^1,\, \dots ,\, \xi^k) \right] \notag \\ &\approx \sum_{p = 1}^{P}\dots \sum_{l = 1}^{L} w_p \dots w_l \, \rho(\xi_p) \dots \rho(\xi_l) \, J_p \dots J_l\,  Q(x,  t; \xi^1_p,\,\dots,\, \xi^k_l) \label{eq:pcm_multi}
\end{align}
where considering each dimension in the random space, we have $k$ summations. The superscript in $\xi^k_l$ indicates the dimension in the random space, and the subscript indicates the collocation point in that dimension. We simplify the notation using $\mathbb{E}\left[ Q(x, t; \xi^1,\, \dots ,\, \xi^k) \right] = \mathbb{E}\left[ Q\right]$, resulting the expression for the standard deviation as

\begin{align}
&\sigma\left[ Q(x, t; \xi^1,\, \dots ,\, \xi^k) \right] = \sigma_{PCM}\left[ Q(x, t; \xi^1,\, \dots ,\, \xi^k) \right] \notag \\&\approx \sqrt{\sum_{p = 1}^{P}\dots \sum_{l = 1}^{L} w_p \dots w_l \, \rho(\xi_p) \dots \rho(\xi_l) \, J_p \dots J_l\, \left(Q(x, t; \xi^1_p,\,\dots,\, \xi^k_l) -  \mathbb{E}\left[ Q\right]\right)^2}. \label{eq:pcm_multi_std}
\end{align}

We assume that the random variables are mutually independent and that the discretization in the parametric space is isotropic. Furthermore, for computational efficiency, a fully tensorial product approach is adequate for this project, as it handles six dimensions or fewer. To avoid the curse of dimensionality in high-dimensional stochastic spaces, Smolyak sparse grids \cite{smolyak1963quadrature} offers a well-regarded solution that reduces the number of realizations while maintaining accuracy. Additionally, methods such as Principal Component Analysis (PCA) \cite{abdi2010principal}, low-rank approximations \cite{chevreuil2015least}, and active subspace methods \cite{constantine2017global} provide effective dimensionality reduction in uncertainty quantification.

\subsection{Sensitivity Analysis}

In our study, we use Sobol indices to assess the global sensitivity of input parameters, following the methodology of \cite{sobol1993sensitivity}. These indices help quantify the contribution of each parameter to the variance of our quantity of interest (QoI). For a detailed derivation of these indices, we refer to the work of \cite{saltelli2010variance}. Let $\xi^j$ represent the $j-$th parameter in our analysis, where $j = 1,2,\dots,k$, and $k$ is the total number of parameters in the space. The impact of $\xi^j$ is measured by its effect on the variance $V$ of the QoI.

\begin{equation}
\label{eq:var}
V_{\xi^j}\left(\mathbb{E}_{\mathbf{\xi}^{\sim j}}(Q | \xi^j)\right)
\end{equation}
where $\mathbf{\xi}^{\sim j}$ represents the set of all parameters except for $\xi^j$, which is held constant at a specific value. This setup involves computing the expected value of the quantity of interest (QoI), $Q$, while fixing $\xi^j$ and subsequently calculating the variance across all possible values of $\xi^j$. According to the Law of Total Variance, we have:

\begin{equation}
\label{eq:law}
V_{\xi^j}\left(\mathbb{E}_{\mathbf{\xi}^{\sim j}}(Q | \xi^j)\right) + \mathbb{E}_{\mathbf{\xi}^{j}}\left(V_{\mathbf{\xi}^{\sim j}}(Q | \xi^j)\right) = V(Q)
\end{equation}

The second term on the left-hand side is the residual term and the total variance is denoted by $V(Q)$. We normalize the equation to compute the first-order sensitivity index which quantifies the proportion of total variance in $Q$ due to the variations in the random variable $\xi^j$. The index is defined as:

\begin{equation}
\label{eq:si}
S_i = \frac{V_{\xi^j}\left(\mathbb{E}_{\mathbf{\xi}^{\sim j}}(U | \xi^j)\right)}{V(U)}
\end{equation}

The sensitivity indices $S_i$ capture the first-order effect of the parameter $\xi^j$ on the variance, excluding interactions between $\xi^j$ and other parameters. Following normalization, the sum of all is less than 1,  $\sum S_i < 1$.  The remainder represents the variance due to higher-order interactions among parameters, which are not addressed in this study but can be analyzed through post-processing of the Polynomial Chaos Method, as discussed by Barros et al. (2021) \cite{barros2021integrated}.

Computing the sensitivity indices $S_i$ can be challenging when uncertainty quantification (UQ) is performed using Monte Carlo (MC) methods due to their computational intensity and need for a large number of samples. However, in this study, we employ the Polynomial Chaos Method as an efficient technique that enables fast and cost-effective computations of global sensitivity. This approach significantly reduces the computational cost associated with traditional MC methods, facilitating more efficient sensitivity analysis.

\subsection{Probability of Failure}

In the final stage of our stochastic analysis, we focus on calculating the failure probability, $p_f$, over time. Typically, methods in reliability literature, such as those discussed by Machado (2015) \cite{machado2015reliability}, rely on Monte Carlo simulations to count failure events to estimate $p_f$. However, stochastic collocation methods present an alternative by providing the moments of a limit state function $g(R, S)$ that must be transformed into a PDF for the computation of $p_f$, a process that involves additional complexity. This transformation can be achieved using the method of moments \cite{low2013new,dang2019novel}, Polynomial Chaos expansions \cite{lasota2015polynomial,garcia2021polinomial}, Gaussian transformations \cite{he2014sparse}, and entropy optimization methods \cite{winterstein2013extremes}. Each of these methods provides a different approach to derive the PDF from the moments, facilitating the computation of $p_f$.

In this study, we introduce an alternative approach that directly estimates the probability of failure, $p_f$, as a standard expectation in uncertainty quantification (UQ). We propose transforming $g$ into a new random variable $h_B$ rather than deriving a probability density function (PDF) from the moments of a limit state function $g$ to determine $P(g < 0)$. The variable $h_B$ is modeled as a Bernoulli random variable with the probability parameter $p_h$. This method simplifies the process by excluding the need for PDF approximation and directly addresses the computation of failure probability.

The definition of $h_B$ comes directly from $g$:

\begin{equation}
\label{eq:h}
h_B = \begin{cases}
0, & \text{if } g \geq 0,\\
1, & \text{otherwise}.
\end{cases}
\end{equation}

In our methodology, each realization of the Polynomial Chaos Method produces a time-series vector $h_B$. Initially, $h_B$ is set to $0$. The value changes to $1$ when $\theta_{\text{max}} > \theta_{\text{lim}}$ or $\varphi_{\text{max}} > \varphi_{\text{lim}}$ and remains at $1$ until the final time-step, effectively forming a step function at the point of failure. At a specific time step, the expectation of $h_B$ reflects the smoothed quantity of interest, with $h_B$ taking on a real value between 0 and 1.

Using the Polynomial Chaos Method to compute the expectation of $h_B$, a Bernoulli random variable, allows us to directly obtain $p_f$ where expectation corresponds to the Bernoulli parameter $p_h$. Thus, a single PCM integration efficiently yields an accurate estimation of $p_f$ at each time step.

\subsection{Numerical Results}
In the deterministic study, we have a clear understanding of how varying initial damage under continuous and unexpected environmental conditions affects the lifespan of transmission lines. Now, we shift our focus to analyzing the impact of parametric uncertainty on maximum temperature or maximum damage in respective scenarios. We model this uncertainty using a uniform distribution with each parameter varying by $10\%$ around its mean value, consistent with the parameters used in the deterministic model.

Our primary goal is to quantify uncertainty and perform sensitivity analysis on parameters affecting maximum damage or maximum temperature. We begin by analyzing the uncertainties in the material parameters  $\xi_m(\omega) = \{A_\sigma(\omega), \gamma(\omega), g_c(\omega), a(\omega)\}$, parameterized due to inaccurate measurement or assumptions in our mathematical modeling. We assume all other material parameters are deterministic. Through UQ and SA, we aim to identify the two parameters within $\xi_m(\omega)$ that have the most significant impact on the variance of $\theta_{\text{max}}$ or $\varphi_{\text{max}}$. After analyzing the material parameters, we examine the uncertainty in loading conditions combined with the two most influential material parameters. These are represented in a separate stochastic space, denoted as $\xi_c(\omega) = {g_c(\omega), a(\omega)\},\theta_b(\omega), w_b(\omega), I_b(\omega), I_A(\omega)}$. From this set, we identify the four most influential parameters based on global SA results. We then combine scenario-specific parameters to understand their influence on the maximum variance.

In all our simulations, we use 5 PCM points per dimension. When analyzing time-series data for maximum temperature across the entire stochastic space, we adjust for varying failure times by truncating the time series at the earliest failure time observed across all realizations.

\subsubsection{Scenario 1 - High Seasonal Winds and Temperature:}

\begin{figure}[H]
	\centering
	\subfloat[Expectation.]{\includegraphics[width=0.45\textwidth]{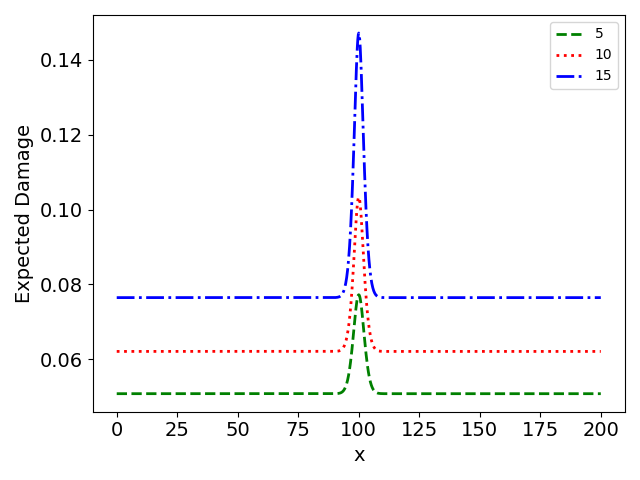}}
	\subfloat[Standard deviation.]{\includegraphics[width=0.45\textwidth]{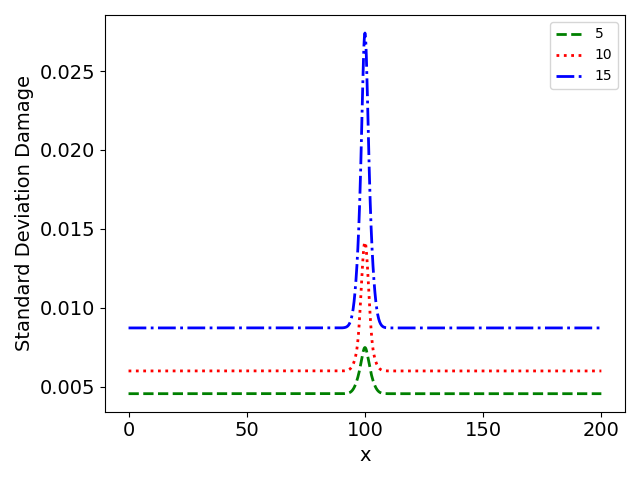}}
	\caption{Expectation and standard deviation damage fields for material parameter uncertainty for set $\xi_m(\omega)$ in 5 years increments (shown in the legend)}
	\label{fig:uq_1_1}
\end{figure}

\begin{figure}[H]
	\subfloat[Expectation.]{\includegraphics[width=0.33\textwidth]{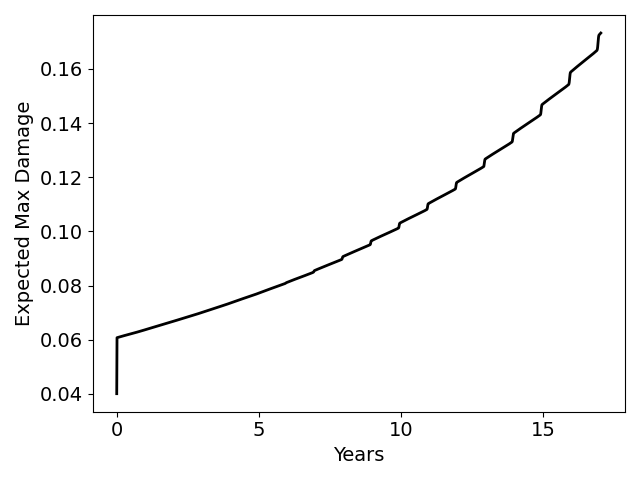}}
	\subfloat[Standard deviation.]{\includegraphics[width=0.33\textwidth]{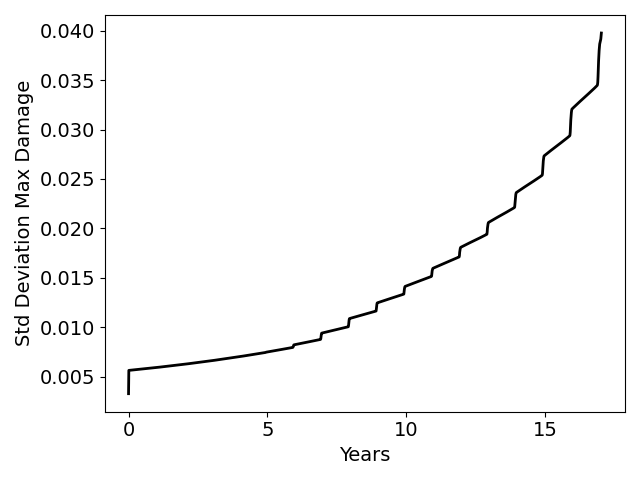}}
	\subfloat[Sensitivity index.]{\includegraphics[width=0.33\textwidth]{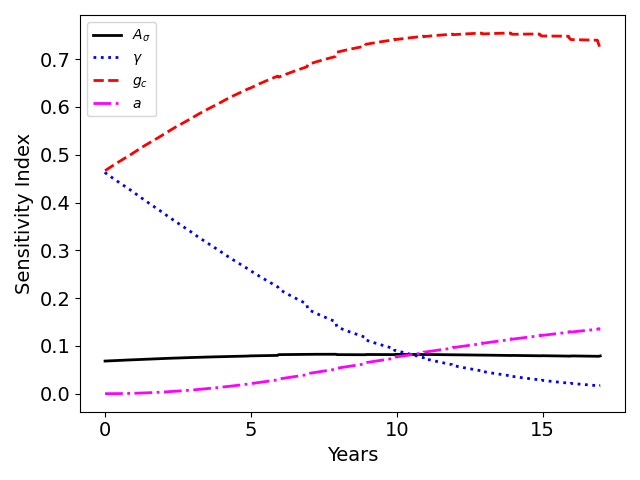}}
	\caption{Time-series of expectation, standard deviation, and sensitivity index of maximum damage under material parametric space set $\xi_m(\omega)$.}
	\label{fig:si_1_1}
\end{figure}

\begin{figure}[H]
	\subfloat[Expectation.]{\includegraphics[width=0.33\textwidth]{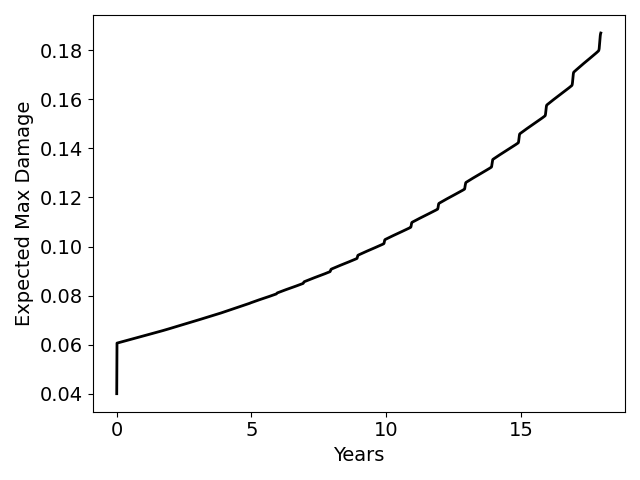}}
	\subfloat[Standard deviation.]{\includegraphics[width=0.33\textwidth]{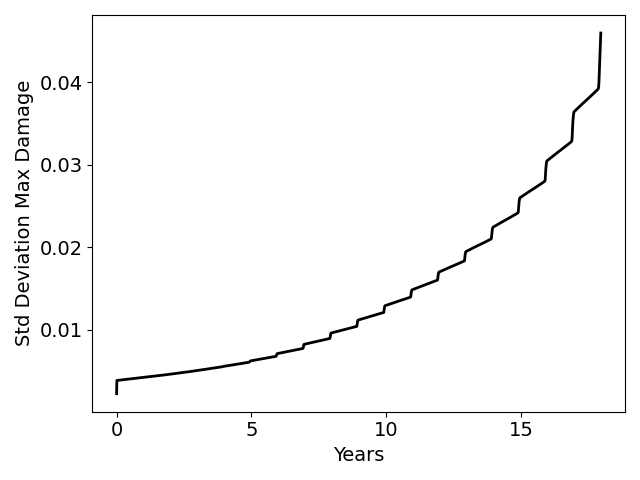}}
	\subfloat[Sensitivity index.]{\includegraphics[width=0.33\textwidth]{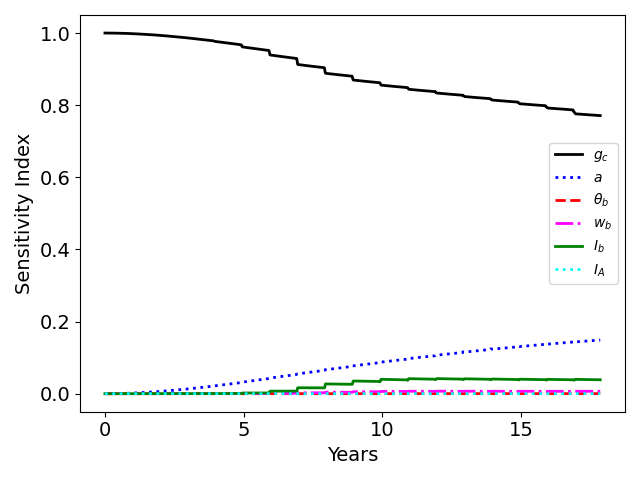}}
	\caption{Time-series of expectation, standard deviation, and sensitivity index of maximum damage under combined parametric space set $\xi_c(\omega)$.}
	\label{fig:si_1_2}
\end{figure}

We start by examining the uncertainty in the material parameters from the set $\xi_m(\omega)$. We analyze the expected damage field and its standard deviation over time, presented in 5-year increments, as shown in Fig.\ref{fig:uq_1_1}. The results confirm that similar to the deterministic solution, the maximum damage and maximum standard deviation occur at the center of the conductor.
Next, we examine the time-series evolution of maximum damage at the center of the transmission line. We then calculate the Sobol indices $S_i$ according to Eq.(\ref{eq:si}) and illustrate the time-series evolution of all parameters from 
$\xi_m(\omega)$ at the center of the transmission line, as shown in Fig.\ref{fig:si_1_1}. Both the expectation and standard deviation are observed to increase over time.
Initially, the damage layer width significantly influences the initiation of damage but becomes less important as the simulation progresses. However, $g_c$ remains the most significant parameter, as it relates to the total energy threshold for fracture. Over time, $a$ becomes increasingly significant due to its relation to the rate of fatigue accumulation. Thus, $g_c$ and $a$ are identified as the two most influential parameters over time. We combine these with the loading parameters to identify the four most influential parameters, as shown in Fig.\ref{fig:si_1_2}. Among the loading conditions, the current is the most influential parameter due to its direct relation to Joule heating. However, over time, material parameters become more significant than loading parameters due to aging effects.

\begin{figure}[H]
	\subfloat[Expectation.]{\includegraphics[width=0.33\textwidth]{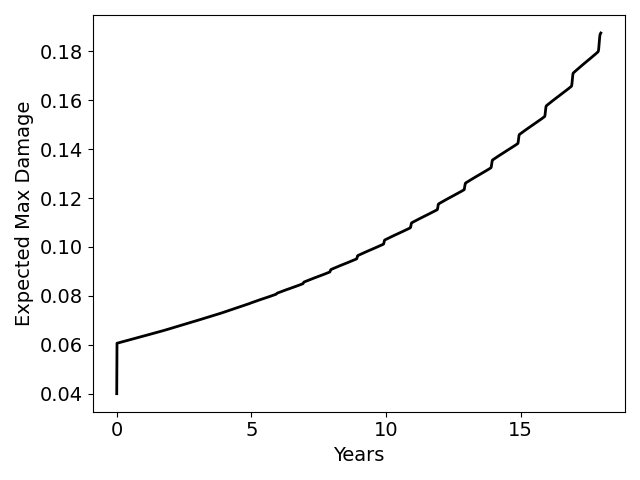}}
	\subfloat[Standard deviation.]{\includegraphics[width=0.33\textwidth]{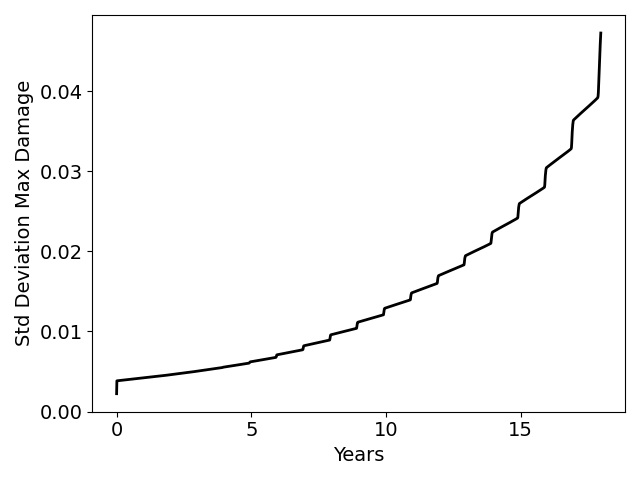}}
	\subfloat[Sensitivity index.]{\includegraphics[width=0.33\textwidth]{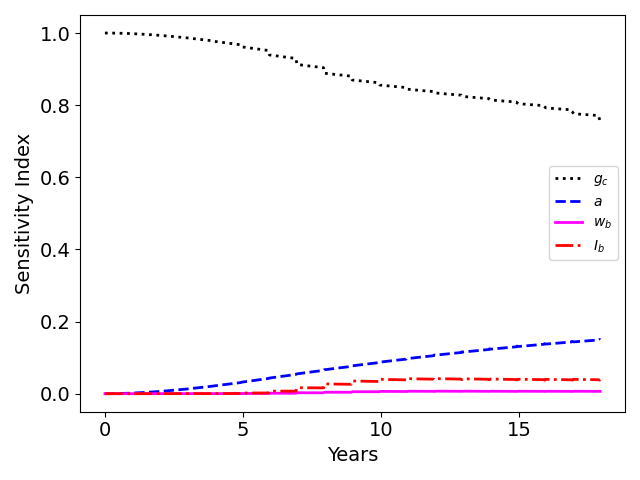}}
	\caption{Time-series of expectation, standard deviation, and sensitivity index of maximum damage under combined parametric space set $\xi_{f1}(\omega)$.}
	\label{fig:si_1_3}
\end{figure}

We then combine the effects of the most influential material and loading condition parameters to form a new set of random parameters $\xi_{f1} (\omega) = \{g_c(\omega), a(\omega),w_b(\omega),I_b(\omega)\}$, and perform a final round of SA. The results are plotted in Fig.\ref{fig:si_1_3}

Finally, we combine the four influential parameters from the final set $\xi_{f1}(\omega)$ with unexpected wind speed, resulting in the set $\xi_1 (\omega) = \{g_c(\omega),a(\omega),w_b(\omega),I_b(\omega),w_{max}\}$. The scenario parameter of high wind speed is crucial when it initially impacts the line, significantly increasing damage. However, its impact reduces after the initial hit as it facilitates convective cooling. Although the extra wind cools the conductor considerably, the expected failure occurs before 10 years.
\begin{figure}[H]
	\subfloat[Expectation.]{\includegraphics[width=0.33\textwidth]{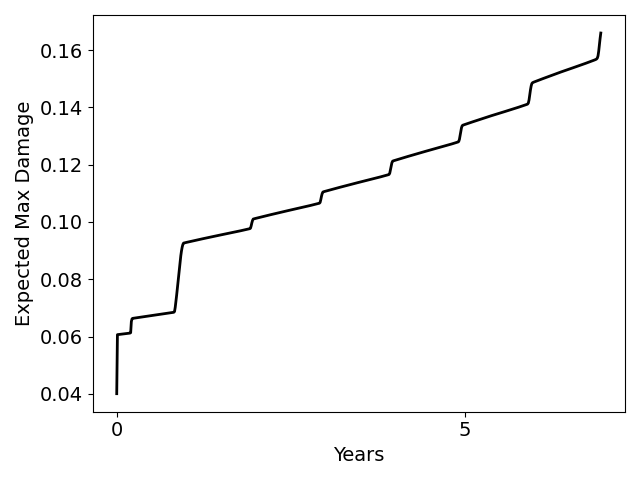}}
	\subfloat[Standard deviation.]{\includegraphics[width=0.33\textwidth]{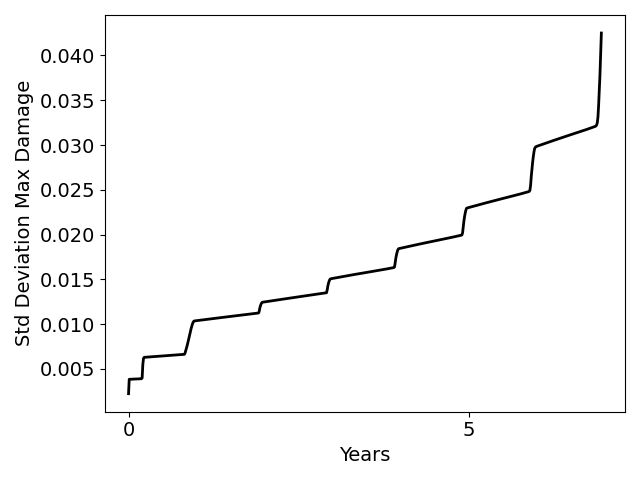}}
	\subfloat[Sensitivity index.]{\includegraphics[width=0.33\textwidth]{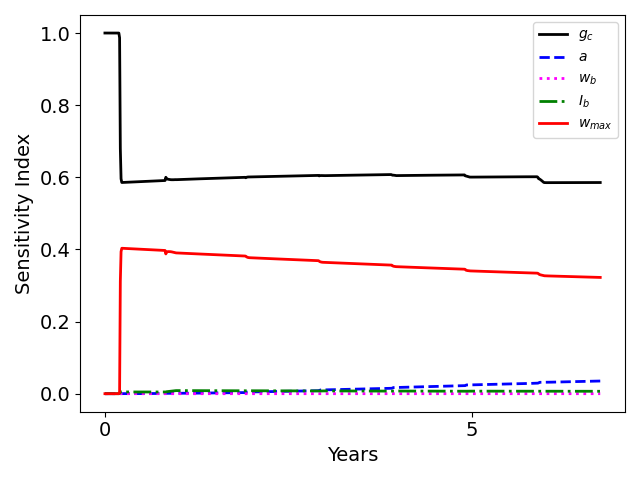}}
	\caption{Time-series of expectation, standard deviation, and sensitivity index of maximum damage under combined parametric space set $\xi_1(\omega)$.}
	\label{fig:si_2_1}
\end{figure}

\subsubsection{Scenario 2 - Wildfire:}
\begin{figure}[H]
	\centering
	\subfloat[Expectation.]{\includegraphics[width=0.45\textwidth]{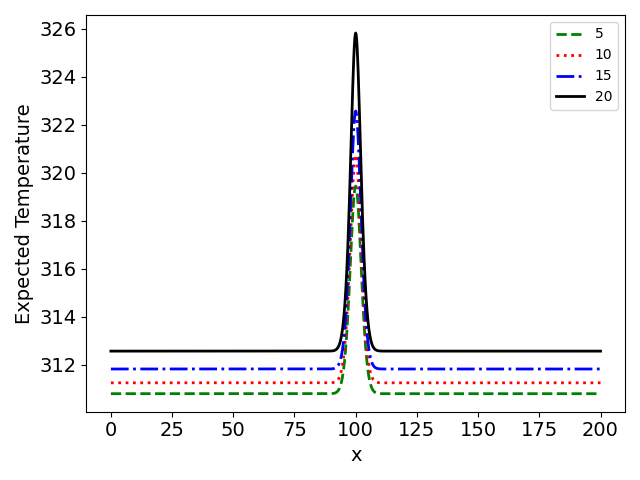}}
	\subfloat[Standard deviation.]{\includegraphics[width=0.45\textwidth]{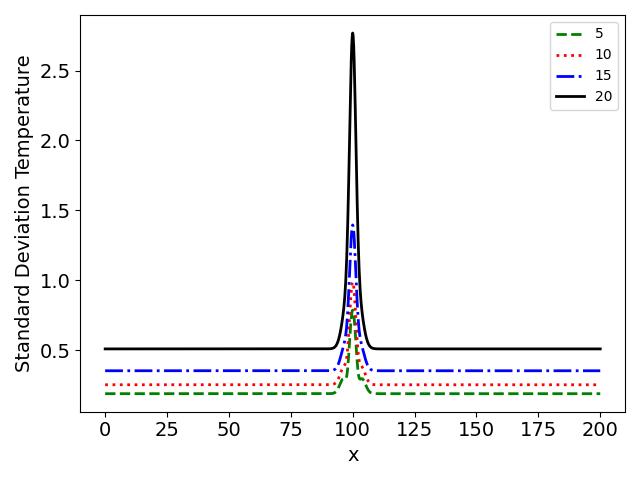}}
	\caption{Expectation and standard deviation temperature fields for material parameter uncertainty for set $\xi_m(\omega)$ in 5 years increments (shown in the legend)}
	\label{fig:uq_3_1}
\end{figure}

We follow a similar procedure to identify the uncertainty in the material parameters from the set $\xi_m(\omega)$ regarding the variance of maximum temperature. The data indicates high temperatures and low wind speeds in the region, leading to the temperature as the mode of failure. We first analyze the expected temperature field and its standard deviation over time, presented in 5-year increments, as shown in Fig.\ref{fig:uq_3_1}. The results reveal that the maximum temperature and maximum standard deviation remain at the center of the conductor. 

\begin{figure}[H]
	\subfloat[Expectation.]{\includegraphics[width=0.33\textwidth]{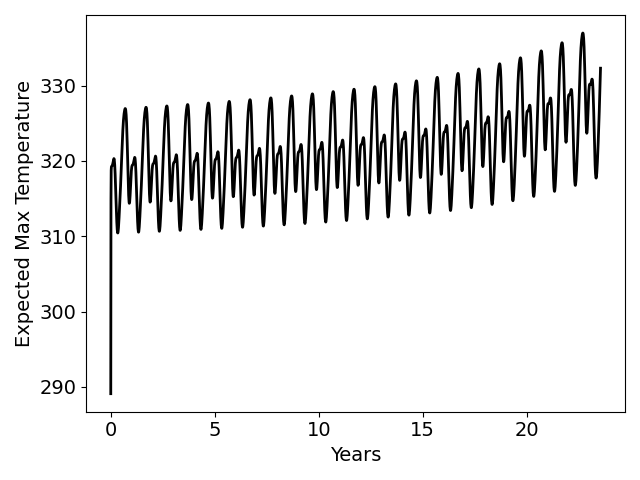}}
	\subfloat[Standard deviation.]{\includegraphics[width=0.33\textwidth]{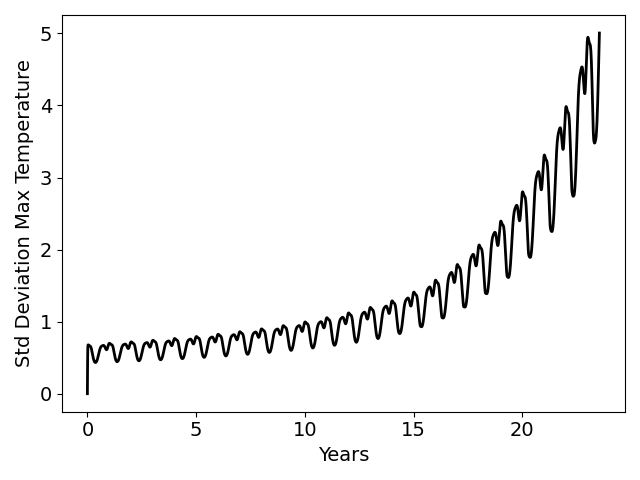}}
	\subfloat[Sensitivity index.]{\includegraphics[width=0.33\textwidth]{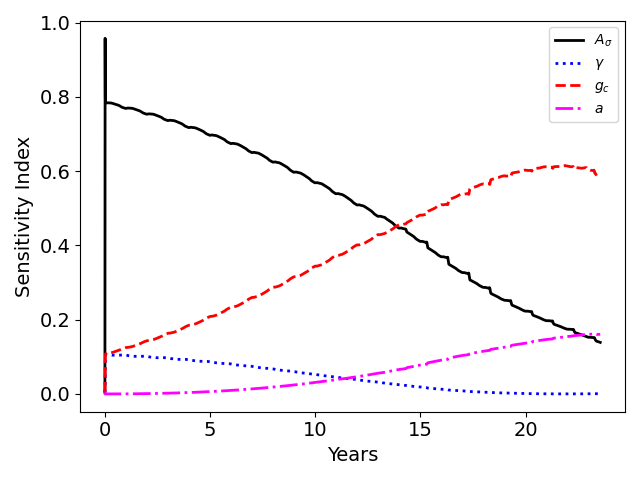}}
	\caption{Time-series of expectation, standard deviation, and sensitivity index of maximum temperature under material parametric space set $\xi_m(\omega)$.}
	\label{fig:si_3_1}
\end{figure}

Next, we examine the time-series evolution of maximum temperature at the center of the transmission line. We then calculate the Sobol indices $S_i$ and illustrate the time-series evolution of all parameters from $\xi_m(\omega)$ at the center of the transmission line, as shown in Fig.\ref{fig:si_3_1}. Both the expectation and standard deviation increase over time, with higher fluctuations.

\begin{figure}[H]
	\subfloat[Expectation.]{\includegraphics[width=0.33\textwidth]{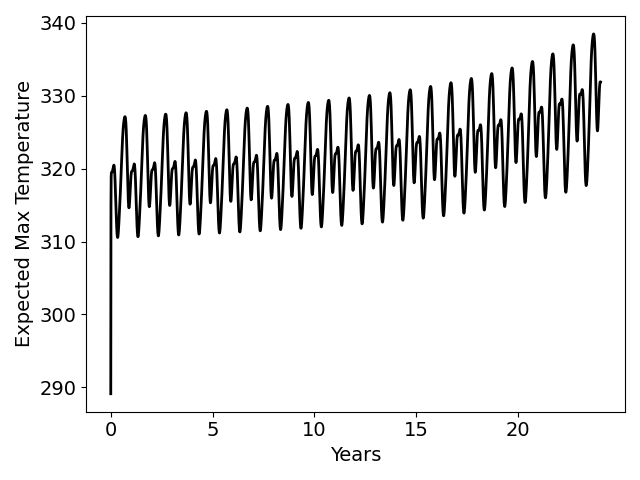}}
	\subfloat[Standard deviation.]{\includegraphics[width=0.33\textwidth]{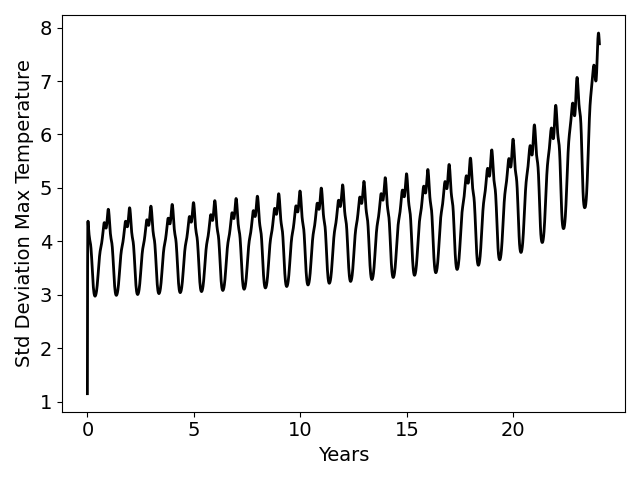}}
	\subfloat[Sensitivity index.]{\includegraphics[width=0.33\textwidth]{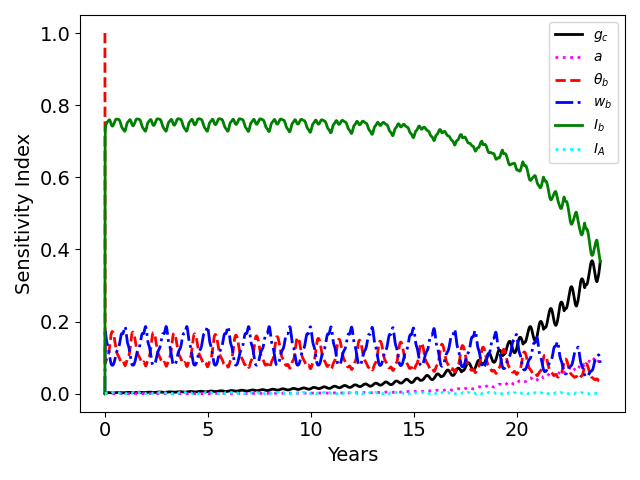}}
	\caption{Time-series of expectation, standard deviation, and sensitivity index of maximum temperature under combined parametric space set $\xi_c(\omega)$.}
	\label{fig:si_3_2}
\end{figure}

Among the material parameters, the cross-sectional area parameter, which drives damage localization, is initially significant in initiating damage but becomes unimportant as the simulation progresses. Over time, $g_c$ and $a$ become more important. We then combine these two most influential material parameters with the loading conditions to identify the four most influential parameters. As shown in Fig.\ref{fig:si_3_2} among the loading parameters, the current base parameter $I_b$ is the most significant in the total uncertainty of $\theta_{max}$, due to its direct correlation with Joule heating. The competition between Joule heating and convective cooling makes the wind base parameter $w_b$ the second most important.

\begin{figure}[H]
	\subfloat[Expectation.]{\includegraphics[width=0.33\textwidth]{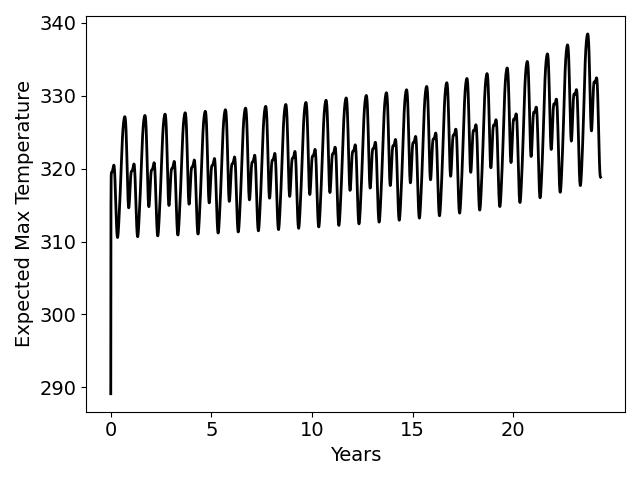}}
	\subfloat[Standard deviation.]{\includegraphics[width=0.33\textwidth]{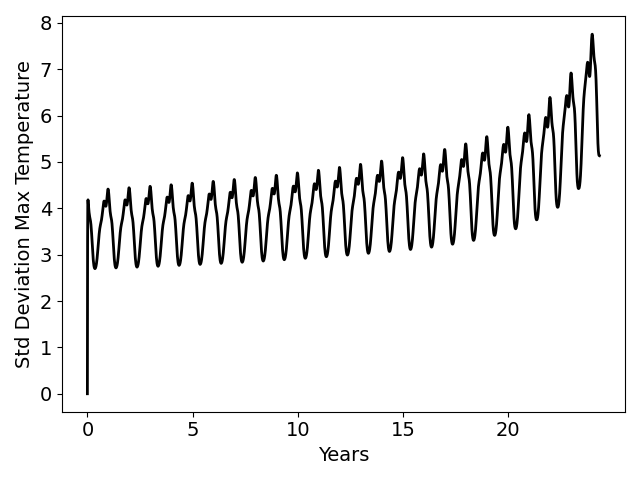}}
	\subfloat[Sensitivity index.]{\includegraphics[width=0.33\textwidth]{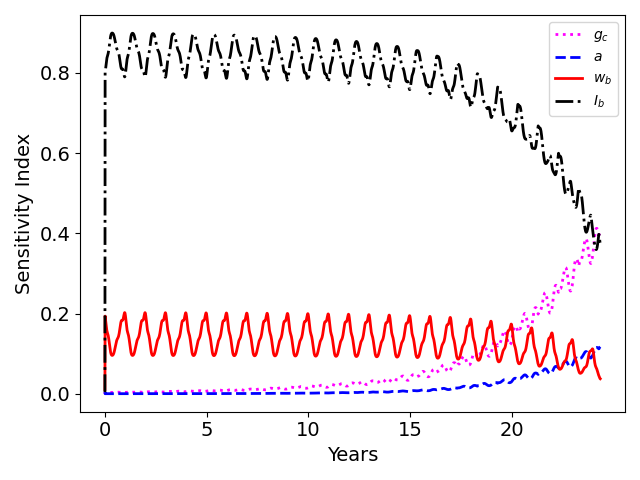}}
	\caption{Time-series of expectation, standard deviation, and sensitivity index of maximum temperature under combined parametric space set $\xi_{f2}(\omega)$.}
	\label{fig:si_3_3}
\end{figure}

\begin{figure}[H]
	\subfloat[Expectation.]{\includegraphics[width=0.33\textwidth]{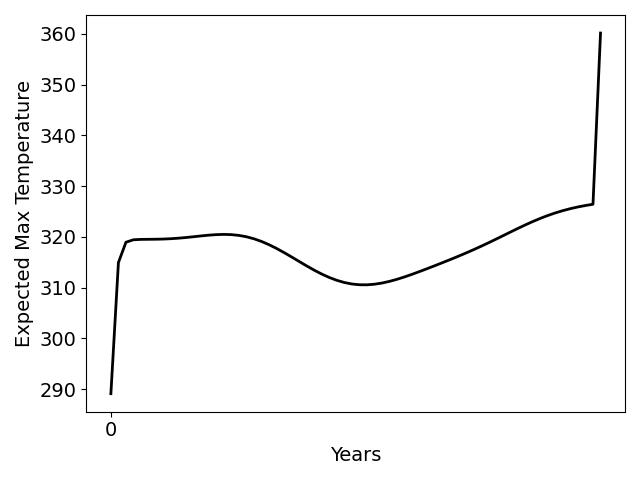}}
	\subfloat[Standard deviation.]{\includegraphics[width=0.33\textwidth]{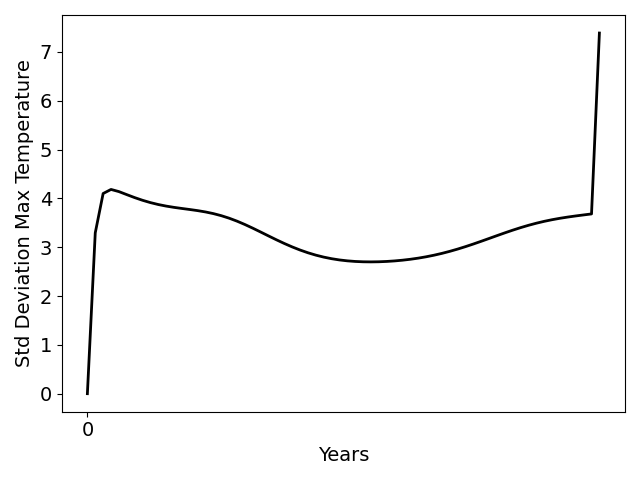}}
	\subfloat[Sensitivity index.]{\includegraphics[width=0.33\textwidth]{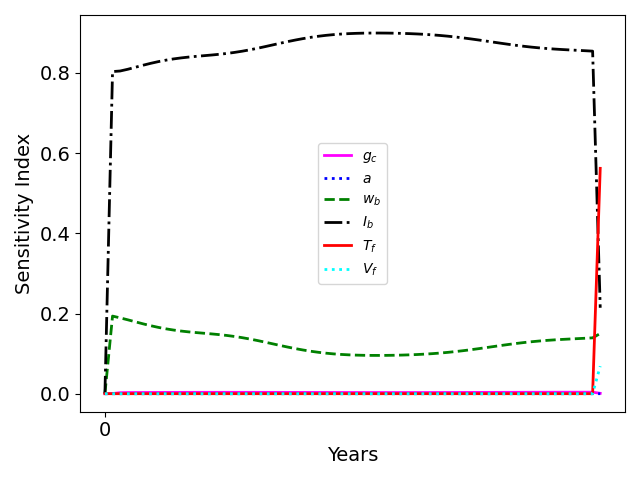}}
	\caption{Time-series of expectation, standard deviation, and sensitivity index of maximum temperature under combined parametric space set $\xi_2(\omega)$.}
	\label{fig:si_4_1}
\end{figure}

Similar to scenario 1, we identify the four influential parameters shown in Fig.\ref{fig:si_3_3} and combine them with two additional parameters: wildfire temperature and view factor. The resulting set is 
$\xi_2 (\omega) = \{g_c(\omega), a(\omega),w_b(\omega),I_b(\omega),T_{fire}, V_f \}$. We then perform the final UQ and sensitivity analysis. As shown in Fig.\ref{fig:si_4_1}, initially, the current is the most significant parameter. However, during the wildfire, the temperature of the wildfire surpasses the current's dominant effect. The maximum temperature and its standard deviation at the center of the transmission line increase suddenly during the wildfire.

\subsubsection{Scenario 3 - Icing:}
\begin{figure}[H]
	\subfloat[Expectation.]{\includegraphics[width=0.33\textwidth]{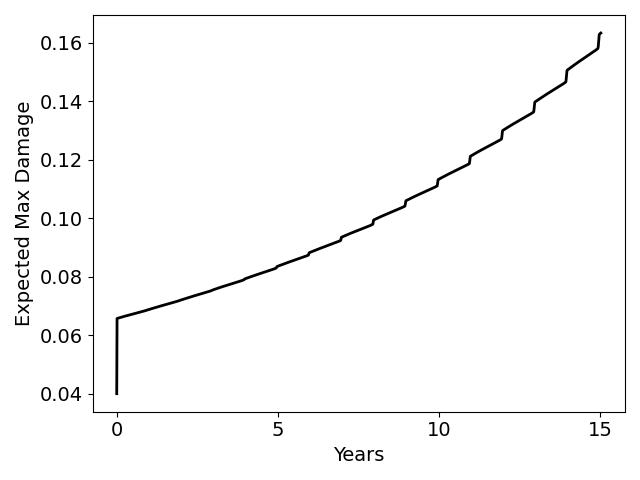}}
	\subfloat[Standard deviation.]{\includegraphics[width=0.33\textwidth]{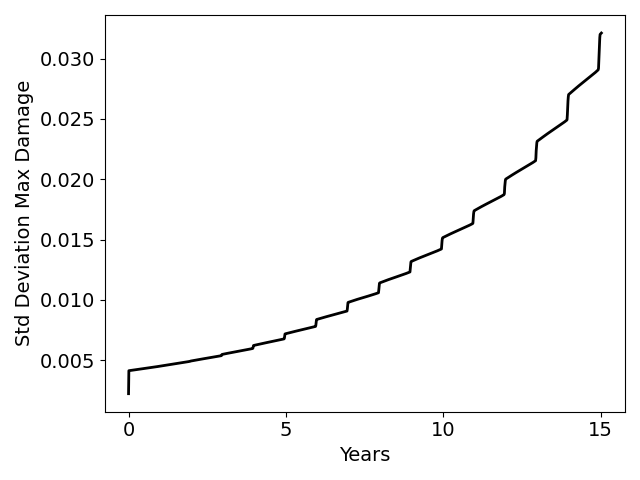}}
	\subfloat[Sensitivity index.]{\includegraphics[width=0.33\textwidth]{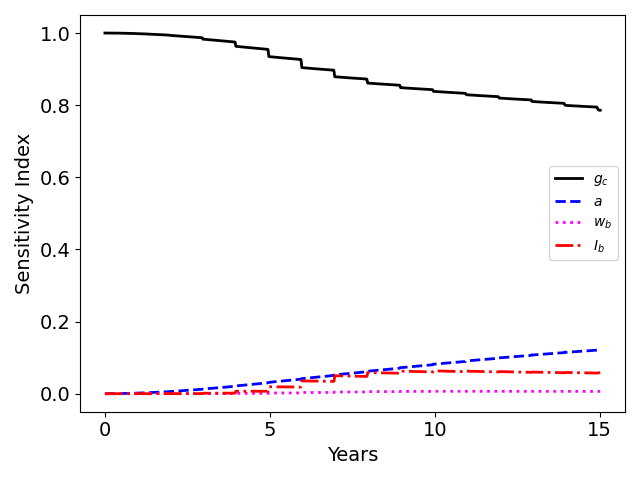}}
	\caption{Time-series of expectation, standard deviation, and sensitivity index of maximum damage under combined parametric space set $\xi_{f3}(\omega)$.}
	\label{fig:si_5_3}
\end{figure}

Similar to scenario 1, the mode of failure is damage. Referencing scenario 1, we combine the most influential parameters for damage to form the final set of model parameters: $\xi_{f3} (\omega) = \{g_c(\omega), a(\omega),w_b(\omega),I_b(\omega)\}$. As shown in Fig.\ref{fig:si_5_3}, the material parameters are the most significant factors accelerating damage in the transmission line. 
We include the thickness of ice as a parameterized loading condition and combine it with the four influential parameters to form the new set 
$\xi_3 (\omega) = \{g_c(\omega), a(\omega),w_b(\omega),I_b(\omega), t_{ice}\}$. We then perform the final UQ and sensitivity analysis to understand the impact of ice on the transmission line. The material parameter $g_c$ remains the most influential. Initially, the thickness of ice has a significant effect due to the mechanical load it imposes on the transmission line. However, over time, this effect diminishes due to the cooling effect.
\begin{figure}[H]
	\subfloat[Expectation.]{\includegraphics[width=0.33\textwidth]{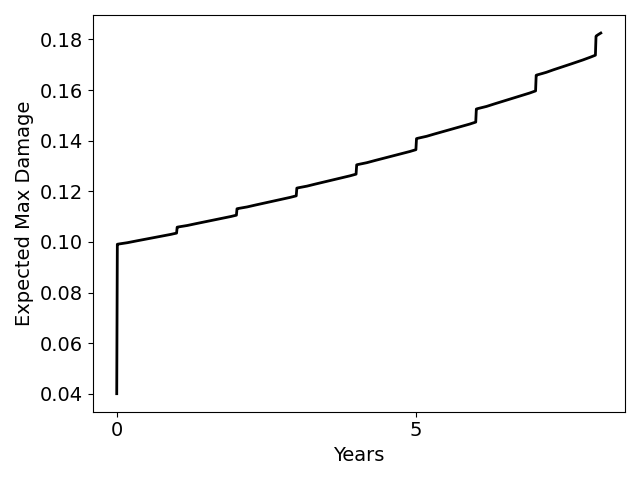}}
	\subfloat[Standard deviation.]{\includegraphics[width=0.33\textwidth]{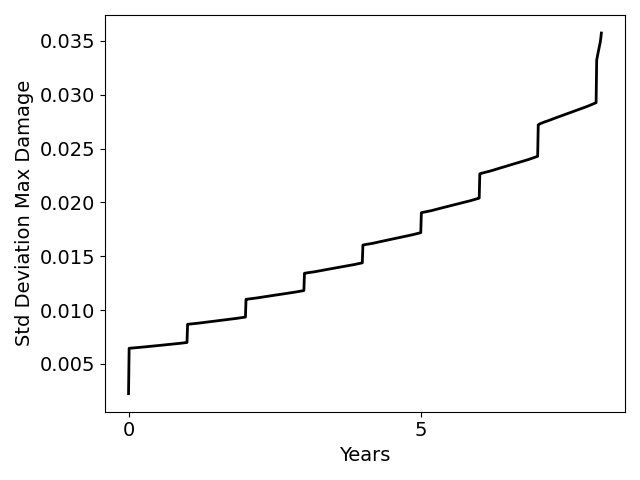}}
	\subfloat[Sensitivity index.]{\includegraphics[width=0.33\textwidth]{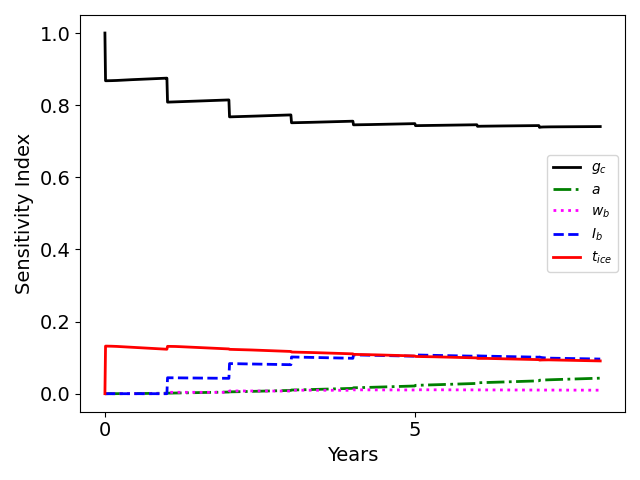}}
	\caption{Time-series of expectation, standard deviation, and sensitivity index of maximum damage under combined parametric space set $\xi_3(\omega)$.}
	\label{fig:si_6_1}
\end{figure}

\subsubsection{Probability of failure}
\begin{figure}[H]
	\centering
	\subfloat[Probability of Failure]{\includegraphics[width=0.45\textwidth]{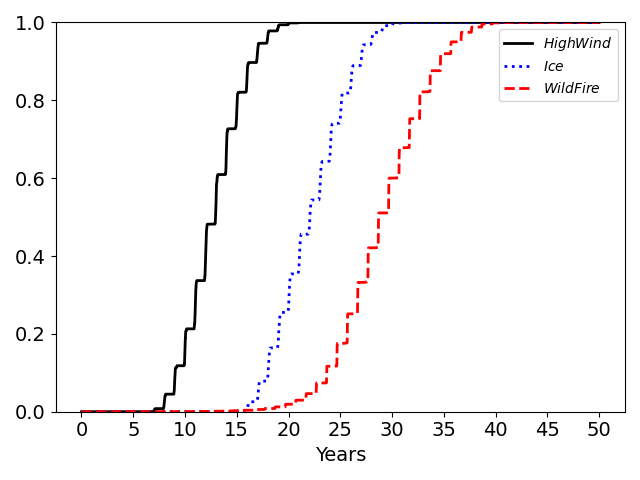}}
	\caption{Probability of failure of three scenarios.}
	\label{fig:pf_all}
\end{figure}

\begin{figure}[H]
	\subfloat[Scenario 1.]{\includegraphics[width=0.33\textwidth]{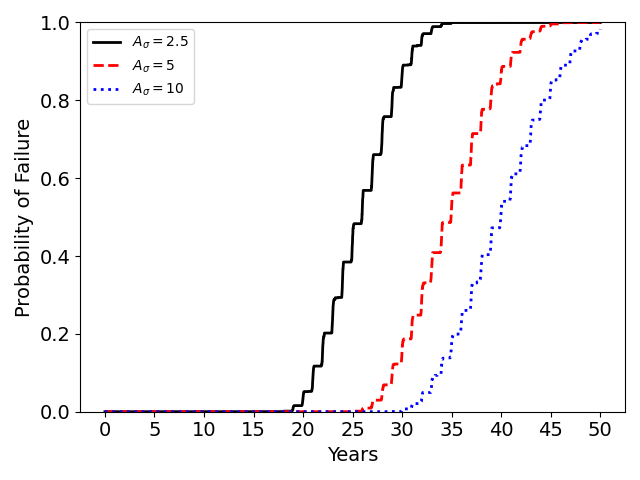}}
	\subfloat[Scenario 2.]{\includegraphics[width=0.33\textwidth]{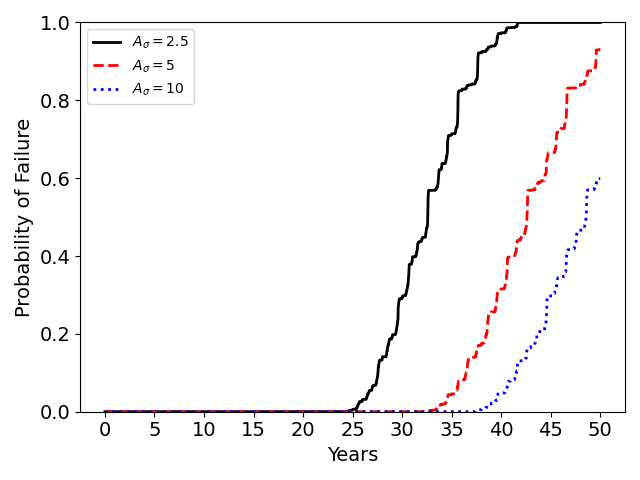}}
	\subfloat[Scenario 3.]{\includegraphics[width=0.33\textwidth]{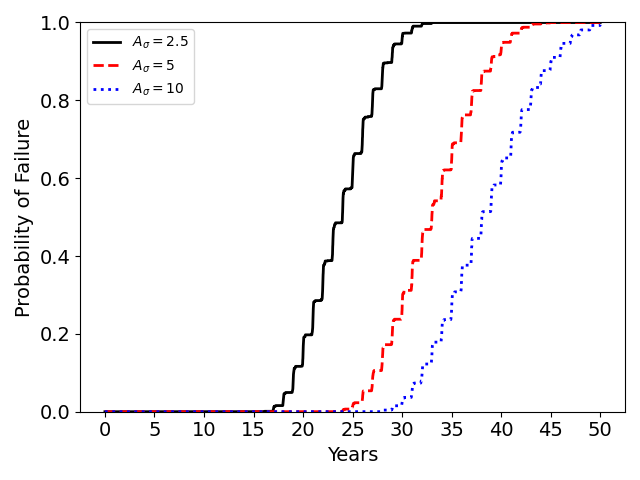}}
	\caption{Probability of failure varying the severity of damage.}
	\label{fig:pf_1}
\end{figure}

\begin{figure}[H]
	\subfloat[Scenario 1.]{\includegraphics[width=0.33\textwidth]{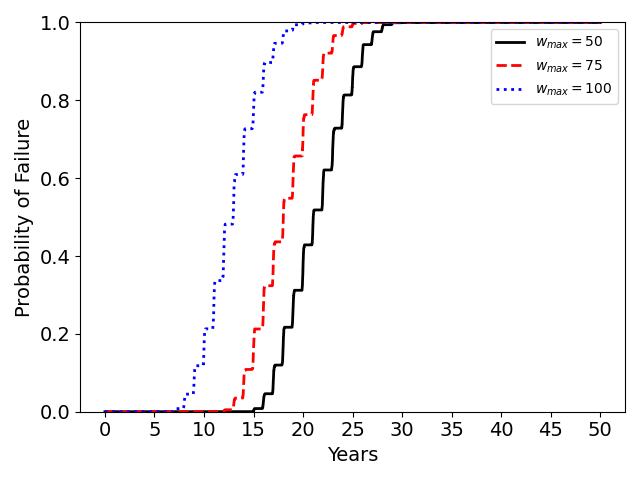}}
	\subfloat[Scenario 2.]{\includegraphics[width=0.33\textwidth]{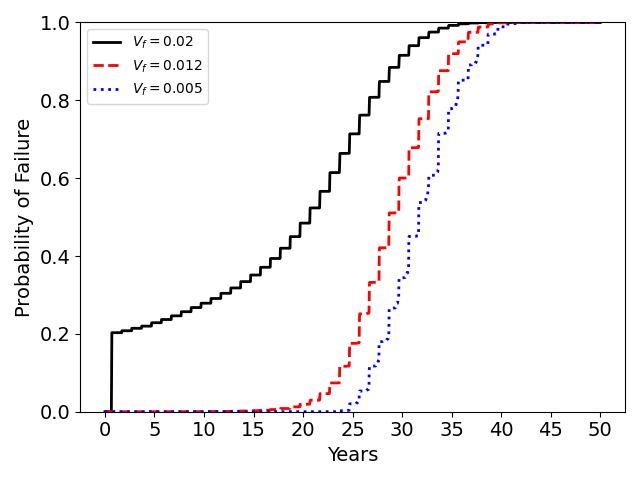}}
	\subfloat[Scenario 3.]{\includegraphics[width=0.33\textwidth]{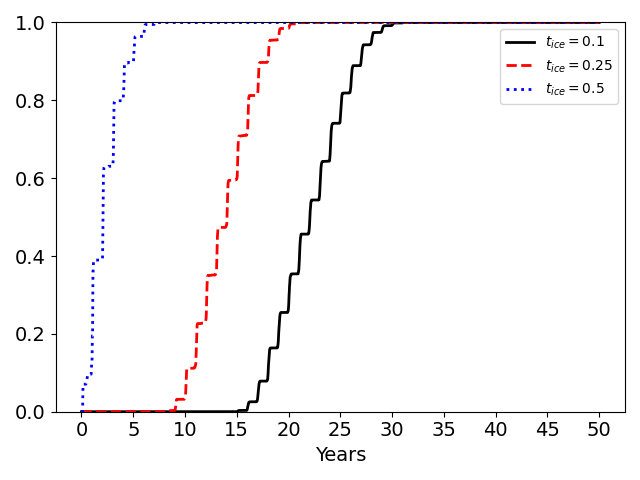}}
	\caption{Probability of failure varying the scenario-specific parameters.}
	\label{fig:pf_2}
\end{figure}

To analyze the probability of failure, we first calculate the expected value of the Bernoulli variable $h_B$ using the PCM with $n=5$ points for each specific scenario. Initially, we compare the probability of failure across three scenarios in Fig.\ref{fig:pf_all}. We use parameter sets $\xi_1(\omega)$, $\xi_2(\omega)$, and $\xi_3(\omega)$, choosing a uniform distribution for the reference mean parameter values with a 10$\%$ lower and upper bound. For these scenarios, we specify a high wind speed of 100 ft/s, an ice thickness of 0.25 inch, and a view factor of 0.0125, assuming severe damage conditions. The probability of failure curve for scenario 1 indicates an early, higher chance of failure. In contrast, scenario 2, influenced by a wildfire, shows a shifted curve to the right, indicating a delayed and reduced occurrence of failure compared to the other scenarios. However, the level of the mean parameters can significantly influence the observed probability of failure. 
To show the impact of mean values, we initially examine the scenarios individually by varying the extent of damage, without taking into account the scenario-specific parameters, as shown in Fig.\ref{fig:pf_1}. The chance of failure increases with increasing severity of damage. For instance, in scenario 2, the presence of severe damage leads to a failure with certainty, whereas less severe damage results in a failure probability of only 50$\%$. 
Finally, we analyze the impact of each scenario-specific parameter for each scenario using sets, $\xi_1 (\omega)$, $\xi_2 (\omega)$, and $\xi_3 (\omega)$. Fig.\ref{fig:pf_2} shows how the mean value of scenario parameters impacts the chance of failure. With the increase in high wind from 50 to 100 ft/s, the chance of failure increases to 90$\%$ at 15 years. Although the high wind speed enhances convective cooling, the impact is more significant leading to the early failure. In scenario 2, the probability of failure in the presence of wildfire increases suddenly if the distance between the transmission line and wildfire is close. The radiative heat transfer is so significant that at a time of wildfire, the probability of failure increases to 20$\%$. We can also infer that if the distance is far enough, the effect of wildfire is insignificant. In the case of scenario 3, the thickness of ice on the transmission line shows a significant reduction in the life span. The failure occurs within 5 years in the presence of heavy ice on the transmission line. Although ice on the layer of lines increases heat transfer through conduction and phase change, the mechanical load it imposes on the conductor is more significant. The combined effect of ice and wind creates a galloping effect reducing the life span of the transmission line.   
\subsubsection{Convergence analysis}
\begin{figure}[H]
	\centering
	\subfloat[PCM.]{\includegraphics[width=0.45\textwidth]{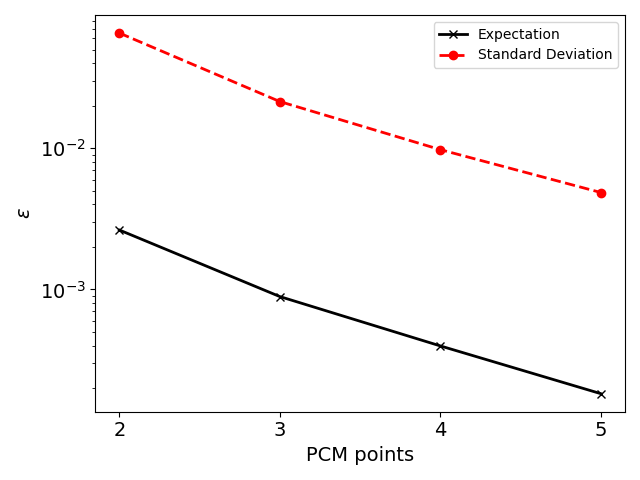}}
      \subfloat[MC.]{\includegraphics[width=0.45\textwidth]{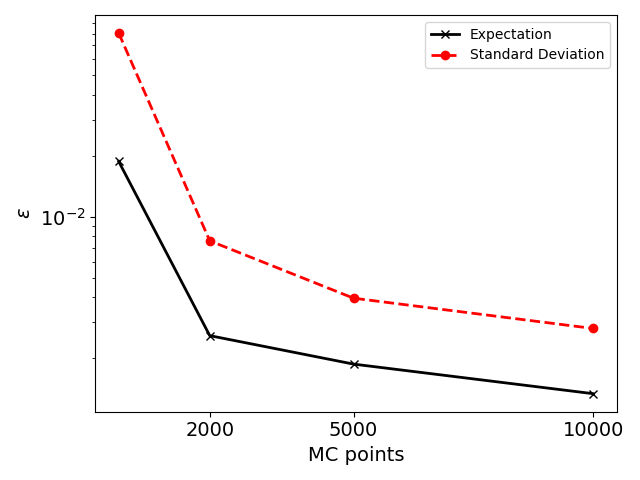}}
	\caption{Convergence Plots.}
	\label{fig:convergence}
\end{figure}

The analysis provides an interesting perspective for understanding the reliability of transmission lines.  However, it is important to verify the consistency of our results. Obtaining an analytical solution would verify the consistency, however, due to the coupling of governing equations, it is not feasible. Therefore, we depend on convergence analysis using the PCM solution as our reference.  
\begin{equation}
\epsilon = \frac{\Vert \varphi -\varphi_{ref}\Vert_2 }{\Vert \varphi_{ref} \Vert_2}.
\end{equation}
We consider scenario 1 as our reference scenario to obtain the convergence plot. Based on our global sensitivity analysis, we identify $g_c$ as the most influential parameter. So we focus on the effects of the most influential parameter $g_c$ simplifying the PCM to a 1-D problem. We consider 100 collocation points as the reference solution. We then compare the PCM solutions to Monte Carlo simulations, analyzing the relative errors in the norm of the damage field at the center of the line at 25 years. The graphical representation is shown in  Fig.\ref{fig:convergence} highlighting the accuracy of PCM relative to 10000 MC realizations.
\section{Conclusion}

We developed a thermo-electro-mechanical model to analyze the reliability of transmission lines in the presence of initial damage and external environmental conditions. We considered historical data on wind and temperature relevant to such scenarios. Our model incorporates a mechanical model for evaluating material damage under prolonged fatigue. It also integrates a thermal model that includes a heat transfer, in terms of Joules heating and convective cooling. However, for scenarios involving ice, the model incorporates heat transfer mechanisms through conduction and phase change. In scenarios with wildfires, it integrates radiative heat transfer. Moreover, the electrical aspect in the model addresses how voltage drops along the line, due to accumulated damage and temperature-induced resistivity. Overall, the model acts as a positive feedback loop, where initially damaged transmission lines deteriorate to the point where the material reaches its threshold limit in terms of either temperature or damage, ultimately leading to failure. 

We used the one-dimensional finite element method to solve a set of governing equations. We considered high wind, wildfire, and icing as three different scenarios to understand their effect on the long-term behavior of the transmission lines. We also varied the damage level to study the effect of the initial damage on the life span of the transmission line. Each scenario was initially analyzed deterministically without considering scenario-specific parameters to understand the effect of damage. Later, the analysis was further extended considering unexpected conditions under varying loading conditions and initial damage. The discrete Fourier transform and Fourier series were used to obtain the cyclic loading condition from the discrete data related to wind and temperature. Current loading was parameterized based on the allowable ampacity to reduce the complexity. Subsequently, we utilized the Probabilistic Collocation Method (PCM) to assess uncertainty quantification (UQ), sensitivity analysis (SA), and probability of failure.

The deterministic solution showed how the scenario-specific condition in the presence of initial damage impacts the failure of the transmission line reducing the life span significantly. The impact of high wind and icing conditions accelerated the accumulation of damage at the center of the line pushing it beyond the threshold leading to early failure. In contrast, the radiative effect of wildfire increased the operating temperature of the line beyond the critical value, failing due to the annealing of the material. 
Under temperature mode of failure, the global sensitivity analysis revealed that electric current, $I_b$ initially has a dominating effect which was later surpassed by fracture energy, $g_c$ over time, and the temperature of wildfire. However, $g_c$ remained the most influencing parameter accelerating the aging although the significance reduced over time. The probability of failure analysis using PCM revealed that the presence of severe initial damage significantly reduced the life of cables. Further, it showed the chance of failure increases in the presence of unexpected conditions such as high wind, wildfire, and icing.  

We acknowledge that our present model incorporates simplifying assumptions that could be further explored in the future. The quasi-static equations could be replaced by transient equations providing a better capture of the evolving nature of the system's response. Furthermore, the one-dimensional simplification of cable could be extended to a catenary shape providing a more realistic scenario. Incorporating current data potentially can improve the prediction precision.  All these changes greatly enhance the reliability evaluations of the transmission line shifting toward a more reliable methodology.

\bibliographystyle{ieeetr}
\bibliography{references}
\end{document}